\documentclass{article}
\usepackage[a4paper]{geometry}
\usepackage[T1]{fontenc}

\usepackage[]{amsmath, amssymb, amsthm, tabularx}

\usepackage[]{a4, xcolor, here}

\usepackage[]{pstricks, pst-text, pst-node, pst-tree, gastex}

\usepackage{tikz}

\usepackage[tikz]{bclogo}

\usepackage{appendix}

\usepackage{boites,graphicx}

\usepackage{soul}

\definecolor{pastelyellow}{rgb}{0.99, 0.99, 0.59}
\definecolor{aqua}{rgb}{0.0, 1.0, 1.0} 
\definecolor{aquamarine}{rgb}{0.5, 1.0, 0.83} 
\definecolor{bananayellow}{rgb}{1.0, 0.88, 0.21}
\definecolor{burgundy}{rgb}{0.5, 0.0, 0.13}
\definecolor{ao(english)}{rgb}{0.0, 0.5, 0.0}

\setul{}{0.2ex}
\setulcolor{bananayellow}

\usepackage[normalem]{ulem}

\usepackage{mathrsfs}

\usepackage{stmaryrd}

\usepackage{enumerate}

\usepackage{paralist}

\usepackage{multienum}

\usepackage{multicol}

\usepackage{fancyhdr}

\usepackage[]{hyperref} 
\hypersetup{
		colorlinks = true,
		linkcolor = red,
		anchorcolor = black,
		citecolor = blue, 
		filecolor = cyan,
		menucolor = red,
		runcolor = cyan,
		urlcolor = blue,
		linkbordercolor = {white},
		linktocpage = true
}

\newtheorem{theorem}{Theorem}[section]
\newtheorem{proposition}[theorem]{Proposition}
\newtheorem{lemma}[theorem]{Lemma}

\theoremstyle{definition}

\newtheorem{example}[theorem]{Example}
\newtheorem{examples}[theorem]{Examples}
\newtheorem{remark}[theorem]{Remark}

\makeatletter
\def\thmhead@plain#1#2#3{%
  \thmname{#1}\thmnumber{\@ifnotempty{#1}{ }\@upn{#2}}%
  \thmnote{ {\the\thm@notefont#3}}}
\let\thmhead\thmhead@plain
\makeatother

\newcommand{\ba}{\mathbf{a}}
\newcommand{\bb}{\mathbf{b}}
\newcommand{\bc}{\mathbf{c}}
\newcommand{\bd}{\mathbf{d}}
\newcommand{\be}{\mathbf{e}}
\newcommand{\bv}{\mathbf{v}}
\newcommand{\bx}{\mathbf{x}}
\newcommand{\by}{\mathbf{y}}

\newcommand{\bG}{\mathbf{G}}
\newcommand{\bH}{\mathbf{H}}

\newcommand{\bS}{\mathbf{S}}
\newcommand{\bT}{\mathbf{T}}

\newcommand{\cC}{\mathcal{C}}
\newcommand{\cD}{\mathcal{D}}

\newcommand{\cF}{\mathcal{F}}
\newcommand{\cG}{\mathcal{G}}
\newcommand{\cH}{\mathcal{H}}

\newcommand{\cS}{\mathcal{S}}
\newcommand{\cU}{\mathcal{U}}
\newcommand{\cV}{\mathcal{V}}

\newcommand{\GL}{\mathrm{GL}}
\newcommand{\SL}{\mathrm{SL}}
\newcommand{\rsp}{\mathrm{rowsp}}
\newcommand{\stab}{\mathrm{Stab}}
\newcommand{\orb}{\mathrm{Orb}}
\newcommand{\rk}{\mathrm{rk}}
\newcommand{\bo}{\mathbf{0}}

\newcommand{\bbF}{{\mathbb F}} 

\renewcommand{\geq}{\geqslant}
\renewcommand{\leq}{\leqslant}

\newcommand{\then}{\Longrightarrow}

\newcommand{\cqd}{\hfill \vrule height 1.5mm width 1.5mm}

{\end{list}}%

\begin{document}

\renewcommand{\headrulewidth}{0pt}

\rhead{ }
\chead{\scriptsize Flag Codes of Maximum Distance and Constructions using Singer Groups}
\lhead{ }

\title{Flag Codes of Maximum Distance\\ and Constructions using Singer Groups}

\author{\renewcommand\thefootnote{\arabic{footnote}}
Miguel \'Angel Navarro-P\'erez\footnotemark[1], 
\renewcommand\thefootnote{\arabic{footnote}} 
 Xaro Soler-Escriv\`a\footnotemark[1]}

\footnotetext[1]{Dpt.\ de Matem\`atiques, Universitat d'Alacant, Sant Vicent del Raspeig, Ap.\ Correus 99, E -- 03080 Alacant. \\\  E-mail adresses: \texttt{miguelangel.np@ua.es, xaro.soler@ua.es.}}

{\small \date{\today}} 

\maketitle

%

\begin{abstract}
In this paper we study flag codes of maximum distance. We characterize these codes in terms of, at most, two relevant constant dimension codes naturally associated to them. We do this first for general flag codes and then particularize to those arising as orbits under the action of arbitrary subgroups of the general linear group. 
We provide two different systematic orbital constructions of flag codes attaining both maximum distance and size. To this end, we use the action of Singer groups and take advantage of the good relation between these groups and Desarguesian spreads, as well as the fact that they act transitively on lines and hyperplanes.
\end{abstract}

\noindent\textbf{Keywords:} Network coding, flag codes, orbit codes, Desarguesian spreads, Singer group actions.


\section{Introduction}\label{sec:Introduction}

{\em Network Coding} appeared in  \cite{AhlsCai00} as a method for maximize the information rate of a network modelled as an acyclic directed multigraph with possibly several senders and receivers. Afterwards, Koetter and Kchischang stated in \cite{KoetKschi08} an algebraic approach for coding in non coherent networks ({\em Random Network Coding}). In this setting, {\em subspace codes} stand as the most appropriate codes for error correction. Given a finite field $\bbF_q$, a subspace code is just a set of subspaces of the $\bbF_q$-vector space $\bbF_q^n$, for a positive integer $n$. Since their definition in \cite{KoetKschi08}, research works on the structure, construction and decoding methods of this type of codes have proliferated considerably, e.g. in \cite{GoManRo2012, GoRa2014,KohnKurz08,ManGorRos2008,TrautRosen18}. In many of these articles, all the codewords are vector subspaces having the same dimension, in which case we speak about {\em constant dimension codes}. Among them, {\em spread codes}, that is, constant dimension codes which partition the ambient space, appear as a remarkable subfamily, since they have the best distance and the largest size for that distance. Besides, a particular way of constructing constant dimension codes is based on the natural and transitive action of the general linear group, $\GL(n,q)$, on the Grassmannian ${\cal G}_q(k,n)$, which is the set of all $k$-dimensional vector subspaces of $\bbF_q^n$ over $\bbF_q$. In this context, an {\em orbit code} is a constant dimension code which is the orbit under the action of some subgroup of $\GL(n,q)$ acting on ${\cal G}_q(k,n)$. When the acting group is cyclic, we call them {\em cyclic orbit codes}. First studied in \cite{TrautManRos2010}, orbit codes have a rich mathematical structure due to the group action point of view. Concerning cyclic orbit codes, {\em Singer groups}, that is, cyclic subgroups of $\GL(n,q)$ of order $q^n-1$, play a very important role. The structure of these groups, as well as their multiple mathematical properties, have permitted to obtain relevant information about  the orbit codes they generate (see \cite{GluMoTro2015,HorTraut2016,ManTrautRosen2011,RoTraut2013,TrautManBraRos2013}, for instance). 
 
{\em Flag codes} can be seen as a generalization of constant dimension codes. In this case, the codewords are flags on $\bbF_q^n$,  that is, tuples of nested vector subspaces of prescribed dimensions. The use of these codes is particularly interesting when there are limitations on the size $q$ of the field or the length $n$ of the information packets to be transmitted  \cite{NobUcho09}. The use of flags in the Network Coding setting began with the work \cite{LiebNebeVaz2018} by Liebhold \textit{et al.} and has continued with several articles that have extended and deepened this line of research \cite{CasoPlanar,CasoGeneral, Kurz20}. The action of the general linear group on the Grassmannian can be easily extended to any variety of flags. Consequently, it also makes sense to study those flag codes arising as orbits of some relevant groups (see  \cite{Cl-MA2020, PlanarOrbConstr, LiebNebeVaz2018}). The paper at hand is also involved in this research. Concretely, we deal with flag codes of maximum distance ({\em optimum distance flag codes}), with special emphasis on those having an orbital structure under the action of a suitable Singer group.

The paper is structured as follows. In Section \ref{sec:Preliminaries} we give all the background we need on finite fields, constant dimension codes and Singer group actions. Section \ref{sec:spread_Singer}  is devoted to highlighting the well-known relationship between Desarguesian spreads and Singer groups, which will be very important for the subsequent derivation of our flag codes. In Section \ref{sec:flagcodes} we start with a summary of known results on flag codes, after which we delve into optimum distance flag codes. In Theorem \ref{teo: characterization ODFC a and b} we characterize these flag codes in terms of at most two constant dimension codes, improving considerably on the, up to now, only known result in this respect (see \cite{CasoPlanar}). Moreover, we give several characterizations of optimum distance flag codes when arising as orbits of the action of an arbitrary subgroup of the general linear group (Theorem \ref{teo:cond_suf_ODFC_orbital}). In Section  \ref{sec:construc}, we use  the action of appropriate Singer groups and the theoretical results previously obtained, in order to provide two systematic orbital constructions of optimum distance flag codes having the best cardinality. First, in Section  \ref{subsect:$n=ks$}, we give a construction starting from a Desarguesian spread of $\bbF_q^n$ and the results set out in Section \ref{sec:spread_Singer}.  Since this kind of construction entails a restriction on the type vector of the flags (see \cite{CasoGeneral}), we complete our research in Section \ref{subsect:$n=2t+1$}, where we deal with the construction of orbit flag codes of full type having maximum distance and size.


\section{Preliminaries}\label{sec:Preliminaries}

\subsection{Finite fields}\label{sec:ff_Singer}

We begin this section by recalling some definitions and results that can be found in any textbook on finite fields (see \cite{Niede}, for instance).

Let $\bbF_q$ denote the finite field of $q$ elements, for a prime  power $q$. Given a positive integer $k$, we put $\bbF_q^{k\times k}$ for the set of all $k\times k$ matrices with entries in $\bbF_q$ and $\GL(k,q)$ for the {\em general linear group} of degree $k$ over $\bbF_q$, composed by all  invertible matrices in $\bbF_q^{k\times k}$. A {\em primitive element}  $\omega$ of the field $\bbF_{q^k}$ is just a generator of the cyclic group $\bbF_{q^k}^*$, which has order $q^k-1$. Let $p(x)=x^k+\sum_{i=0}^{k-1}p_ix^i\in\bbF_q[x]$ be the minimal polynomial of $\omega$ over $\bbF_q$. It turns out that $p(x)$ is the characteristic polynomial of the matrix

\[
M_k=\begin{pmatrix}
 0 & 1        &  0        &  \cdots & 0\\
  0      & 0  &  1        &  \cdots & 0\\
  \vdots &  \vdots  &   \vdots  &  \ddots & \vdots \\
   0      & 0        &  0  &  \cdots & 1\\
  -p_0      & -p_1        &  -p_2        &  \cdots & -p_{k-1}
\end{pmatrix}\in \GL(k,q),
\] 
which is called the {\em companion matrix} of $p(x)$. In particular, $M_k$ can be seen as a root of $p(x)$ and the finite field $\bbF_{q^k}$ can be realized as $\bbF_{q^k}\cong \bbF_{q}[\omega]\cong\bbF_q[M_k]$, 
where the last field isomorphism is given by:
\begin{equation}\label{eq:def_phi}
\begin{array}{cccc}
\phi : & \bbF_{q}[\omega]                 &  \longrightarrow &         \bbF_{q}\left[ M_k \right]\\
         & \sum_{i=0}^{k-1}a_i\omega^i       & \longmapsto & \sum_{i=0}^{k-1}a_iM_k^i ,
\end{array}\end{equation}
Consequently, the multiplicative order of $M_k$ is $q^k-1$, that is, $M_k$ generates a cyclic subgroup $\langle M_k\rangle$ of order $q^k-1$ in  $\GL(k,q)$. Equivalently, $M_k$ is a primitive element of the finite field $\bbF_q[M_k] \subseteq \bbF_{q}^{k\times k}$.

For any positive integer $n$ and $i\in\{0,\dots, n\}$, we denote by ${\cal G}_{q}(i,n)$ the set of all $i$-dimensional vector subspaces of $\bbF_q^n$ over $\bbF_q$, which is called the {\em Grassmann variety} (or simply the {\em Grassmannian}). For any  positive integer $s\geq 2$, the field isomorphism $\phi$ is useful to map vector subspaces of $\bbF_{q^k}^s$ into vector subspaces of $\bbF_q^{ks}$.  For $m\in\{1,\dots, s\}$, one has the following embedding map 
\begin{equation}\label{eq:def_varphi}
\begin{array}{cccc}
\varphi : & {\cal G}_{q^k}(m,s)                &  \longrightarrow &         {\cal G}_{q}(km,ks)\\
         &\rsp \begin{pmatrix}
 a_{11}    &  \cdots & a_{1s}\\
  \vdots &   \ddots  & \vdots \\
  a_{m1}  &    \cdots & a_{ms}
             \end{pmatrix}       
                                              & \longmapsto &\rsp \begin{pmatrix}
                                              \begin{array}{c|c|c}
                                              \phi(a_{11}) &   \cdots & \phi(a_{1s})\\ \hline
 									          \vdots &  \ddots & \vdots \\ \hline
  									          \phi(a_{m1})  &    \cdots & \phi(a_{ms})
                                              \end{array}
             								\end{pmatrix}, 
\end{array}
\end{equation}
which is called a {\em field reduction map}. In particular, $\varphi(\cal U)\cap \varphi(\cal V)=\varphi(\cal U\cap\cal V)$, for all ${\cal U}, {\cal V}\in {\cal G}_{q^k}(m,s)$, since $\varphi$ is an injective map.

Besides,  we can also use $\phi$ to obtain the following group monomorphism (see \cite[Th. 2.4]{Geertrui2016}): 
\begin{equation}\label{eq:def_psi}
\begin{array}{cccc}
\psi : & \GL(s,q^k)                &  \longrightarrow &         \GL(ks,q)\\
         & \begin{pmatrix}
 a_{11}    &  \cdots & a_{1s}\\
  \vdots &   \ddots  & \vdots \\
  a_{s1}  &    \cdots & a_{ss}
             \end{pmatrix}       
                                              & \longmapsto &  \begin{pmatrix}
                                              \begin{array}{c|c|c}
                                              \phi(a_{11}) &   \cdots & \phi(a_{1s})\\ \hline
 									          \vdots &  \ddots & \vdots \\ \hline
  									          \phi(a_{s1})  &    \cdots & \phi(a_{ss})
                                              \end{array}
             								\end{pmatrix}.
\end{array}
\end{equation}

\subsection{Subspace codes, orbit codes and Singer group actions}\label{subsec:subspaces}

For any integers $1\leq k < n$, the Grassmannian $\cG_q(k, n)$ can be seen as a metric space endowed with the following  \emph{subspace distance} (see \cite{KoetKschi08}):

\begin{equation}\label{eq:dist}
d_S(\cU,\cV)=\dim(\cU+\cV)-\dim(\cU\cap \cV)=2(k -\dim(\cU\cap \cV)),
\end{equation}
for all $\cU, \cV\in\cG_q(k, n)$. Using this metric, a \emph{constant dimension code} is just a nonempty subset $\cC$ of $\cG_q(k, n)$ and its \emph{minimum distance} is defined as 
\begin{equation}\label{eq:distC}
d_S(\cC)= \min\{ d_S(\cU, \cV) \ | \ \cU, \cV\in \cC, \ \cU\neq \cV\}\leq \left\{\begin{array}{lll}
                                                                                                                 2k  & \mbox{if} & 2k\leq n,\\
                                                                                                                 2(n-k) & \mbox{if} & 2k\geq n.
                                                                                                                 \end{array}\right.
\end{equation}
In case that $|\cC|=1$, we put $d_S(\cC)=0$. In any other case, $d_S(\cC)$ is always a positive even integer. When the upper bound provided by (\ref{eq:distC}) is attained we say that $\cC$ is a {\em constant dimension code of  maximum distance}. 

Notice that the distance $d_S(\cC)=2k$ can only be attained if $2k\leq n$ and different subspaces in $\cC$ pairwise intersect trivially. This class of codes of maximum distance were introduced in \cite{GoRa2014} as {\em partial spread codes} since they generalize the class of {\em spread codes}, previously studied in \cite{ManGorRos2008}. A spread code  is  just a spread in the geometrical sense, that is, its elements pairwise intersect trivially and cover the whole space $\bbF_q^n$ (see \cite{Hirschfeld98}). The size of a partial spread code of dimension $k$ (or $k$-partial spread code) is always upper bounded by 
\begin{equation}\label{eq:cota_partial_spread}
\frac{q^n-q^r}{q^k-1},
\end{equation}
where $r$ is the reminder obtained dividing $n$ by $k$ (see \cite{GoRa2014}). In turn, $k$-spreads exist if, and only if, $k$ divides $n$ and, in this case, they have cardinality $(q^n-1)/(q^k-1)$, which is the largest size among $k$-partial spreads of $\bbF_q^n$.  

On the other hand,  if $2k\geq n$ and $\cC\subseteq \cG_q(k, n)$, then the {\em dual code} of $\cC$ is the set $\cC^\perp=\{{\cal V}^\perp\ |\ {\cal V}\in \cC\}.$ It is a constant dimension code of dimension $n-k$ with  the same cardinality and distance than $\cC$ (see \cite{KoetKschi08}). In particular, if $d_S(\cC)=2(n-k)$, then $\cC^\perp$ is an $(n-k)$-partial spread code and the size of $\cC$ can be also upper bounded in terms of (\ref{eq:cota_partial_spread}). 

Notice that any code $\cC$ included in the Grassmannian of lines ${\cal G}_q(1,n)$ or in the Grassmannian of hyperplanes ${\cal G}_q(n-1,n)$, with $|\cC|\geq 2$,  is a constant dimension code of maximum distance and its size is upper bounded by $(q^n-1)/(q-1)$.

{\em Orbit codes} were introduced in \cite{TrautManRos2010} as constant dimension codes arising as orbits under the action of some subgroup of the general linear group. Given a subspace $\cV\in\cG_q(k, n)$ and a full-rank matrix $V\in \bbF_q^{k\times n}$ generating $\cV$, that is, $\cV=\rsp(V)$, the map
\[
\begin{array}{ccc}
{\cal G}_q(k,n)  \times \GL(n,q) &  \longrightarrow  & {\cal G}_q(k,n)\\
(\cV, A)                                      & \longmapsto       & \cV\cdot A=\rsp(VA),
\end{array}
\]
is independent from the choice of $V$ and it defines a group action on $\cG_q(k, n)$ (see \cite{TrautManRos2010}).
For a subgroup $\bH$ of $\GL(n,q)$, the {\em orbit code} $\orb_{\bH}(\cV)$ is just:
\[
\orb_{\bH}(\cV)=\{ \cV\cdot A\ |\ A\in\bH\}\subseteq   {\cal G}_q(k,n).
\]
The size of an orbit code can be computed as 
$
|\orb_{\bH}(\cV)|=\frac{|\bH|}{|\stab_{\bH}(\cV)|},
$
where $\stab_{\bH}(\cV)=\{A\in\bH\ |\  \cV\cdot A=\cV\}$ is the stabilizer subgroup of the subspace $\cV$ under the action of $\bH$. If $\bH=\stab_{\bH}(\cV)$, then $\orb_{\bH}(\cV)=\{\cV\}$ and $d_S(\orb_{\bH}(\cV))=0$. In any other case, the minimum distance of the orbit code $\orb_{\bH}(\cV)$ can  be computed as (see  \cite{TrautManRos2010}) : 
\[
d_S(\orb_{\bH}(\cV))=\min\{d_S(\cV, \cV\cdot A)\ |\  A\in\bH\setminus \stab_{\bH}(\cV) \}.
\] 
When the acting group is cyclic, the corresponding orbit codes are called {\em cyclic orbit codes}. These particular type of orbit codes have been deeply studied in several papers (see \cite{GluMoTro2015, ManTrautRosen2011, RoTraut2013,TrautManBraRos2013} for instance). 

In this paper we will use the action of {\em Singer cyclic subgroups} of $\GL(n,q)$, which are generated by the so called {\em Singer cycles} of $\GL(n,q)$. These are elements of $\GL(n,q)$ having order $q^n-1$, which is the largest element order in $\GL(n,q)$. Although Singer cycles are not necessarily conjugate, the generated subgroups are always conjugate subgroups of $\GL(n,q)$ (see \cite{Hes,HupI} for more information on Singer groups).

The following seminal result about the action of a Singer group on the Grassmannian of lines and the Grassmannian of hyperplanes is due to Singer (1938) and will be used extensively throughout this article:

 \begin{theorem}\cite[Th. 6.2]{BethJungLenz1999}\label{teo:stabSinger_rectes} 
Any Singer cyclic subgroup $\bS$ of $\GL(n, q)$ acts transitively on both $\cG_{q}(1,n)$ and $\cG_{q}(n-1,n)$. Moreover, for any $l\in \cG_{q}(1,n)$ and any $h\in\cG_{q}(n-1,n)$, it holds
\[
\stab_{\bS}(l)=\stab_{\bS}(h)=\{aI_s\ |\ a\in\bbF_{q}^*\},
\]
which is the unique cyclic subgroup of $\bS$ of order $q-1$.
\end{theorem}


 \section{Desarguesian spread codes and Singer groups}\label{sec:spread_Singer}
 
 In this section we focus on the action of Singer groups in order to obtain certain $k$-spreads of $\bbF_q^n$ as their orbits, for $n=ks$ and $s\geq 2$. Later on, in Section \ref{subsect:$n=ks$}, we will use the results obtained here to construct orbit flag codes by considering Singer groups and their subgroups. 

Consider the field reduction map $\varphi$ defined in $(\ref{eq:def_varphi})$ which maps vector subspaces of $\bbF_{q^k}^s$ into vector subspaces of $\bbF_{q}^n$. Applying $\varphi$ to a constant dimension code $\cC\subseteq \cG_{q^k}(m,s)$, we obtain another constant dimension code $\varphi(\cC)\subseteq \cG_{q}(km,n)$.  Since $\varphi$ is injective, it preserves intersections and therefore it follows that $d_{S}(\varphi(\cC))=kd_{S}(\cC)$. In particular, if $\cC$ attains the maximum possible distance, then  $\varphi(\cC)$ do it as well. Even more, if $\cC$ is a $m$-spread code of $\bbF_{q^{k}}^{s}$, then $\varphi(\cC)$ is a $km$-spread code of $\bbF_{q}^{n}$.

We will use two constant dimension codes of $\bbF_{q}^{n}$ constructed in this way. First, from the spread of lines of $\bbF_{q^{k}}^{s}$, we consider 
\begin{equation}\label{def:spreadS}
\cS=\varphi(\cG_{q^k}(1,s))\subseteq \cG_{q}(k,n),
\end{equation}
which is a $k$-spread of $\bbF_q^n$. Originally due to Segre (see \cite{Segre64}), in the network coding setting, this construction appears for the first time in \cite{ManGorRos2008}.  Secondly, from the Grassmannian of hyperplanes of $\bbF_{q^{k}}^{s}$, we obtain 
\begin{equation}\label{def:codi_n-k}
\cH=\varphi(\cG_{q^{k}}(s-1,s))\subseteq \cG_{q}(n-k,n),
\end{equation}
which is a constant dimension code of  $\bbF_q^n$ with maximum distance. 

Notice that the field reduction map $\varphi$ together with the group monomorphism $\psi$ defined  in $(\ref{eq:def_psi})$ make it possible to establish the following relation between  the group action of $\GL(s,q^k)$ on $\cG_{q^k}(m,s)$ and the group action of $\GL(n,q)$ on $\cG_q(km,n)$:

\begin{equation}\label{eq:equiv_accions}
\varphi(\cV\cdot A)=\varphi(\cV)\cdot \psi(A),
\end{equation}
for all $\cV\in\cG_{q^k}(m,s)$ and $A\in \GL(s,q^k)$. In particular, we will use this equality  to relate the respective actions of two Singer cyclic subgroups in which we are very interested. 

Let $\alpha$ be a primitive element of $\bbF_{q^{n}}$. Recall that, given the minimal polynomial of $\alpha$ over $\bbF_{q^k}$ and its companion matrix, $M_s\in \GL(s,q^k)$, then $\bbF_{q^{n}}\cong \bbF_{q^k}[\alpha]\cong\bbF_{q^k}[M_s]$. 
Therefore, the multiplicative order of $M_s$ is $q^{n}-1$ and  $\bbF_{q^k}\left[M_s\right]=\{0_{s\times s}\} \cup \langle M_s\rangle\subseteq \bbF_{q^{k}}^{s\times s}$. In particular, $M_s$ is a Singer cycle of $\GL(s,q^k)$, generating the Singer cyclic subgroup $\langle M_s\rangle$ of $\GL(s,q^k)$.
Furthermore,  $\psi(\langle M_s \rangle)=\langle \psi(M_s)\rangle$ is a Singer cyclic subgroup of $\GL(n,q)$. These two Singer groups will be crucial in Section \ref{subsect:$n=ks$}. 

Coming back to the spread $\cS=\varphi(\cG_{q^k}(1,s))$ defined in ($\ref{def:spreadS})$, let us denote $\cG_{q^k}(1,s)=\{l_{1}, l_{2},\dots, l_{r}\}$ and  $\cS=\{\cS_{1},\cS_{2},\dots, \cS_r\}$, where $r=\frac{q^n-1}{q^k-1}$ and  $\cS_{i}=\varphi(l_{i})$ for all $i=1,\dots, r$.  In accordance with (\ref{eq:equiv_accions}) and Theorem \ref{teo:stabSinger_rectes}, for every $i\in\{1,\dots, r\}$, we can write  
\begin{eqnarray}\label{eq:S_orbSinger}
\cS & = & \varphi(\cG_{q^k}(1,s))=\varphi(\orb_{\langle M_s\rangle }(l_{i}))=\{\varphi(l_{i}\cdot A)\ |\ A\in\langle M_s\rangle \}\nonumber\\
      & = & \{\cS_{i} \cdot \psi(A) \ |\ \psi(A)\in\langle \psi(M_s)\rangle\}= \orb_{\langle \psi(M_s)\rangle}(\cS_{i}).
\end{eqnarray}
In an analogous way, we denote $\cG_{q^k}(s-1,s)=\{h_{1}, h_{2},\dots ,h_{r}\}$ and the code $\cH=\{\cH_{1},\cH_{2},\dots, \cH_r\}$, where $\cH_{i}=\varphi(h_{i})$, for all $i=1,\dots, r$. Then, for every $h_i\in \cG_{q^k}(s-1,s)$, we obtain that 
\begin{equation}\label{eq:H_orbSinger}
\cH=\varphi(\cG_{q^k}(s-1,s))=\varphi(\orb_{\langle M_s\rangle}(h_{i}))=\orb_{\langle \psi(M_s)\rangle }(\cH_{i}).
\end{equation}
That is, the transitive action of $\langle M_s\rangle $ on the lines and hyperplanes of $\bbF_{q^k}^s$ is translated to the transitive action of $\langle \psi(M_s)\rangle $ on the constant dimension codes of maximum distance  $\cS$ and $\cH$. Moreover, from Theorem \ref{teo:stabSinger_rectes} we also obtain that, for any $\cS_{i}\in\cS$ and $\cH_{i}\in\cH$ it holds

\begin{eqnarray}\label{eq:stab_G}
\stab_{\langle \psi(M_s)\rangle}(\cS_{i}) & = & \stab_{\langle \psi(M_s)\rangle}(\cH_{i})  =  \psi\left( \stab_{\langle M_s \rangle}(l_i) \right)  =  \{\psi(aI_s)\ |\ a\in\bbF_{q^k}^*\} \nonumber\\
  & = & \left\{  \begin{pmatrix}
 													\begin{array}{c|c|c}
 													\phi(a) &    \cdots & 0_{k\times k}\\ \hline
 													\vdots  &    \ddots & \vdots \\ \hline
  														0_{k\times k}   &    \cdots & \phi(a)
  													 \end{array}
             								\end{pmatrix} \ \Big|\ a\in\bbF_{q^k}^*
 \right\},
\end{eqnarray}
which has order $q^k-1$ (see also \cite[Lemma 3.5]{Geertrui2016}).

\begin{remark}\label{rem: desarguesian spread}
Following \cite{LavrauwGeertrui2018} and \cite[Cor. 3.8]{Geertrui2016}, we can say that a $k$-spread ${\cD}$ of $\bbF_q^n$ is {\em Desarguesian} if  it is $\GL(n,q)$-equivalent to $\cS$, that is, if there exists $B\in \GL(n,q)$ such that 
\[
\cD=\cS\cdot B=\{\cS_{1},\cS_{2},\dots,  \cS_r\}\cdot B=\{\cS_{1}\cdot B,\cS_{2}\cdot B,\dots,  \cS_r\cdot B\}.
\]
From \cite{Drudge2002,Geertrui2016}, we know that, given a Singer cyclic subgroup of $\GL(n,q)$, there exists a unique $k$-spread of $\bbF_q^n$ which appears as its orbit. Moreover this $k$-spread is Desarguesian. Consequently,  $\cS$ is the unique  $k$-spread of $\bbF_{q}^{n}$ which arises as an orbit of the Singer cyclic subgroup $\langle \psi(M_s)\rangle$. Besides, given another Desarguesian $k$-spread $\cD$ of $\bbF_{q}^{n}$ and $B\in \GL(n,q)$ such that $\cS\cdot B=\cD$, it follows that
\begin{eqnarray}\label{eq:spread_des_orbital}
\cD    & = & \cS \cdot B = \orb_{\langle \psi(M_s)\rangle}(\cS_{i}) \cdot B \nonumber\\
       & = & \{ \cS_{i} \cdot A\ |\ A\in\langle \psi(M_s)\rangle\}\cdot B \nonumber\\
       & = & \{ \cS_{i} \cdot  AB\ |\ A\in \langle \psi(M_s)\rangle\}\nonumber\\
       & = & \{ (\cS_{i} \cdot B) \cdot B^{-1} AB\ |\ A\in \langle \psi(M_s)\rangle\}\nonumber\\
       & = & \orb_{\langle \psi(M_s)\rangle^{B}}(\cS_{i} \cdot B),
\end{eqnarray}
that is, the $k$-spread $\cD$ appears as the orbit under the action of the Singer cyclic subgroup $\langle \psi(M_s)\rangle^{B}=B^{-1}\langle \psi(M_s)\rangle B$. 
In Section \ref{subsect:$n=ks$}, we will construct flag codes with maximum distance and an orbital structure. To do so, we will make use of the $k$-spread $\cS$ defined in ($\ref{def:spreadS})$ and its orbital structure under the action of the Singer group $\langle \psi(M_s)\rangle$. By virtue of (\ref{eq:spread_des_orbital}), to work with any other Desarguesian $k$-spread $\cD$, it would be enough to consider the group $\langle \psi(M_s)\rangle^{B}$ instead of $\langle \psi(M_s)\rangle$. 
\end{remark}


\section{On Flag codes}\label{sec:flagcodes}
The present section is devoted to a theoretical study of flag codes. We start in Section \ref{subsec:prel_flags} with a revision of some known results. Next, in Section \ref{sec:flags} we focus on flag codes attaining the maximum possible distance and give a characterization of them that considerably improves the one obtained in \cite{CasoPlanar}. We finish the section by studying how to construct these flag codes as orbits (or union of orbits) of arbitrary subgroups of the general linear group.

\subsection{Background of flag codes}\label{subsec:prel_flags}
Given integers $0<t_1<\dots <t_r<n$, a {\em flag} $\mathcal{F}=(\mathcal{F}_1,\dots,  \mathcal{F}_r)$ of type $(t_1, \dots, t_r)$ on $\bbF_q^n$ is a sequence of nested subspaces of $\bbF_q^n$, 
\[
\{0\}\subsetneq \mathcal{F}_1 \subsetneq \cdots \subsetneq \mathcal{F}_r \subsetneq \mathbb{F}_q^n,
\]
with $\mathcal{F}_i \in \mathcal{G}_q(t_i,n)$, for all $i=1,\dots, r.$ We say that $\mathcal{F}_i$ is the {\em $i$-th subspace} of the flag $\cF$ and when the type vector is $(1, 2, \dots, n-1)$ we speak about {\em full flags}. 

The set of all flags of type $(t_1, \dots, t_r)$ on $\mathbb{F}_q^n$ is known as the  {\em flag variety} of type $(t_1, \dots, t_r)$ on $\mathbb{F}_q^n$  and will be denoted here by $\mathcal{F}_q((t_1,\dots, t_r),n)$.

The subspace distance defined in (\ref{eq:dist}) for the Grassmann variety can be naturally extended to  $\mathcal{F}_q((t_1,\dots, t_r),n)$ as follows. Given  $\cF=(\mathcal{F}_1,\dots,  \mathcal{F}_r)$ and $\cF'=(\mathcal{F}'_1,\dots,  \mathcal{F}'_r)$ two flags in $\mathcal{F}_q( (t_1, \dots, t_r),n)$, their {\em flag distance} is
\[
d_f(\cF,\cF')= \sum_{i=1}^r d_S(\mathcal{F}_i, \mathcal{F}'_i).
\] 
The use of flags in the network coding setting appears for the first time in \cite{LiebNebeVaz2018}. Since then, several papers on this subject have recently appeared (see, for instance, \cite{Cl-MA2020,PlanarOrbConstr,CasoPlanar, CasoGeneral, Kurz20}). If $\emptyset\neq \cC\subseteq \mathcal{F}_q((t_1,\dots, t_r),n)$, we say that $\cC$ is a {\em flag code of type $(t_1,\dots,t_r)$} on $\bbF_q^n$. The {\em minimum distance} of $\cC$ is given by
\[
d_f(\cC)=\min\{d_f(\cF,\cF')\ |\ \cF,\cF'\in\cC, \ \cF\neq \cF'\}.
\]
As for subspace codes, if $|\cC|=1$, we put $d_f(\cC)=0$. Notice that $d_f(\cC)$ is upper bounded by (see \cite{CasoPlanar})
\begin{equation}\label{eq:quotamaxdistflag}
d_f(\mathcal{C}) \leq 2\left(\sum_{2t_i \leq n} t_i + \sum_{2t_i > n} (n-t_i)\right) .
\end{equation}
For any $i\in\{1,\dots,r\}$, the  $i$-{\em projected code} $\cC_i$ of $\cC$ is defined in \cite{CasoPlanar} as the constant dimension code
\[
\cC_i=\{ {\cF}_i \ |\   (\mathcal{F}_1,\dots,  \mathcal{F}_r) \in \cC \}\subseteq \cG_q(t_i,n).
\]
Observe that $|\cC_i|\leq |\cC|$, for every $1\leq i\leq r$. When we have the equalities $|\cC|=|\cC_1|=\dots=|\cC_r|$, the flag code $\cC$ is said to be {\em disjoint} (see \cite{CasoPlanar}).
In the same paper, {\em optimum distance flag codes} (ODFCs, for short) are defined as flag codes attaining the upper bound given in $(\ref{eq:quotamaxdistflag})$. The close relationship between a flag code $\cC$ and its projected codes can be used to characterize this class of flag codes by using the concept of disjointness.

\begin{theorem} \cite[Th. 3.11]{CasoPlanar}\label{teo:carac_odfc}
Let $\cC$ be a flag code of type $(t_1, \ldots, t_r)$. The following statements are equivalent:
\begin{enumerate}[$(i)$]
\item $\cC$ is an ODFC.
\item  $\cC$ is a disjoint flag code and every projected code $\cC_i$ is a constant dimension code of maximum distance.
\end{enumerate}
\end{theorem}

As a result, the size of an ODFC can be upper bounded in terms of the upper bound given by (\ref{eq:cota_partial_spread}) of Section \ref{subsec:subspaces}.

\begin{theorem}\cite[Th. 3.12]{CasoPlanar}\label{teo:ODFC_maxcard}
Let $k$ be a divisor of $n$ and assume that $\cC$ is an ODFC of type $(t_1,\dots, t_r)$ on $\bbF_q^n$. If $k$ is a dimension in the type vector, say $t_i=k$, then $|\cC|\leq \frac{q^n-1}{q^k-1}$. Equality holds if, and only if, the projected code $\cC_i$ is a $k$-spread of $\bbF_q^n$. 
\end{theorem}

On the other hand, requiring an ODFC to have a $k$-spread as one of its projected codes leads to a condition on its type vector.

\begin{theorem}\cite[Th. 3.3]{CasoGeneral}\label{teo:max_ad_type}
Let $\cC$ be an ODFC of type $(t_1,\dots, t_r)$ on $\bbF_q^n$. Assume that some dimension $t_i=k$ divides $n$ and the associated projected code $\cC_i$ is a $k$-spread. Then, for each $j\in\{1,\dots,r\}$, either $t_j\leq k$ or $t_j\geq n-k$ holds.
\end{theorem}

Consequently, when $n=ks$ and $s\geq 2$, the {\em full admissible type vector} for an ODFC having a $k$-spread as a projected code is $(1,\ldots, k,n-k,\ldots,n-1)$. We will come back to this situation in Section \ref{subsect:$n=ks$}. Besides, notice that for $s=2$ the full admissible type vector is just the full type vector. A construction of ODFCs with the largest possible size in this particular case, together with a decoding algorithm, can be found in \cite{CasoPlanar}.

The action of the general linear group on the Grassmannian seen in Section \ref{sec:Preliminaries} can be naturally extended to flags as follows. This approach already appears in \cite{Cl-MA2020, PlanarOrbConstr, LiebNebeVaz2018}.
Given a flag $\cF=(\mathcal{F}_1,\ldots,  \mathcal{F}_r)\in \mathcal{F}_q((t_1,\ldots, t_r),n)$ and a subgroup $\bH$ of $\GL(n,q)$, the {\em orbit flag code} generated by $\cF$ under the action of $\bH$ is
\[
\orb_{\bH}(\cF)=\{ \cF\cdot A\ |\ A\in\bH\}=\{(\mathcal{F}_1\cdot A,\ldots,  \mathcal{F}_r\cdot A)\ |\ A\in\bH\} \subseteq   \mathcal{F}_q((t_1,\ldots, t_r),n).
\]
Its associated stabilizer subgroup is $\stab_{\bH}(\cF)=\{A\in \bH\ |\  \cF\cdot A=\cF\}$ and it holds $|\orb_{\bH}(\cF)|=\frac{|\bH|}{|\stab_{\bH}(\cF)|}$. Moreover, the minimum distance of the orbit flag code can be obtained as 
\[
d_f(\orb_{\bH}(\cF))=\min\{d_f(\cF, \cF\cdot A)\ |\  A\in\bH\setminus \stab_{\bH}(\cF) \}
\] 
and it holds $d_f(\orb_{\bH}(\cF))=0$ if, and only if, $\stab_\bH(\cF)=\bH$. The projected codes of an orbit flag code are orbit (subspace) codes. More precisely, for every $1\leq i\leq r$, we have $\orb_{\bH}(\cF)_i=\orb_{\bH}(\cF_i)\subseteq \cG_q(t_i,n).$ Besides, the stabilizer subgroup of $\cF$ is closely related to the ones of its subspaces:

\begin{equation}\label{eq:inter_stab}
\stab_{\bH}(\cF)=\bigcap_{i=1}^r \stab_{\bH}(\cF_i).
\end{equation}
Remark that, fixed an acting subgroup $\bH$ of $\GL(n,q)$, the cardinality of the flag code $\orb_{\bH}(\cF)$ and their projected codes are determined by the orders of the corresponding stabilizer subgroups of $\bH$. The equality given in $(\ref{eq:inter_stab})$ allows to obtain that disjoint flag codes, in the orbital scenario, involve an equality of subgroups and not only of cardinalities.

\begin{proposition}\cite[Prop. 3.5]{PlanarOrbConstr}\label{prop:disjorb}
Given a flag $\cF=(\cF_1, \ldots, \cF_r)$ of type $(t_1, \ldots, t_r)$ on $\bbF_q^n$ and a
subgroup $\bH$ of $\GL(n,q)$, the following statements are equivalent:
\begin{enumerate}[$(i)$]
    \item $\orb_\bH(\cF)$ is a disjoint flag code. 
    \item $\stab_\bH(\cF)=\stab_\bH(\cF_1)=\cdots= \stab_\bH(\cF_r)$. 
    \item $\stab_\bH(\cF_1)=\cdots= \stab_\bH(\cF_r)$. 
\end{enumerate}
\end{proposition}


\subsection{On Optimum Distance Flag Codes}\label{sec:flags}

Next, we will go deeper into the theoretical study of ODFCs, obtaining some important properties that will be useful for the subsequent orbital constructions (Section \ref{sec:construc}). In Theorem \ref{teo:carac_odfc}, ODFCs are characterized in terms of all their projected codes. In this section, we go one step further and present a new criterion to characterize them just regarding, at most, two of its projected codes. To achieve this result, we start studying how the fact of having a projected code with maximum distance gives us some information about the cardinality and distance of some other projected codes.

\begin{proposition}\label{prop: a projected code with max dist}
Let $\cC$ be a flag code of type $(t_1, \dots, t_r)$ on $\bbF_q^n$ having a projected code $\cC_i$ of maximum distance, for some $i\in\{1,\dots,r\}$, and take flags $\cF, \cF'\in\cC$ such that $\cF_{i}\neq\cF'_{i}$.
\begin{enumerate}[$(a)$]
\item If \ $2t_i\leq n$, then $d_S(\cF_j, \cF'_j)=2t_j$, for every $1\leq j\leq i.$ In particular, $|\cC_i|\leq |\cC_{j}|$, for  $1\leq j\leq i$.
\item If\  $2t_i \geq n$, then $d_S(\cF_j, \cF'_j)=2(n-t_j)$ for all $i\leq j\leq r.$ As a consequence,  $|\cC_i|\leq |\cC_{j}|,$ for values $i\leq j\leq r$.
\end{enumerate}
\end{proposition}
\begin{proof}
Assume that $\cC_i$ is a constant dimension code of maximum distance. In particular, $|\cC_i|> 1$ and we can consider flags $\cF, \cF'\in\cC$ such that $\cF_{i}\neq\cF'_{i}$. We distinguish two possibilities:
\begin{enumerate}[$(a)$]
\item If $2t_i\leq n$, then $d_S(\cF_i, \cF'_i)=d_S(\cC_i)=2t_i$ or, equivalently, $\cF_i\cap\cF'_i=\{0\}$. Hence, for every $1\leq j\leq i$, it holds that $\cF_j\cap\cF'_j=\{0\}$ and then $d_S(\cF_j, \cF'_j)= 2t_j$, which is the maximum possible distance between $t_j$-dimensional subspaces. In particular, $\cF_j\neq\cF'_j$ and we conclude that $|\cC_i|\leq |\cC_j|$, for all $1\leq j\leq i$.
\item If $2t_i\geq n$, then $d_S(\cF_i, \cF'_i)=d_S(\cC_i)=2(n-t_i)$ or, equivalently, $\cF_i+\cF'_i=\bbF_q^n$. As a result, when we consider higher dimensions $t_j\geq t_i$ in the type vector, we obtain $\cF_j+\cF'_j=\bbF_q^n$ as well. Consequently, $d_S(\cF_j, \cF'_j)=2(n-t_j)$, which is the maximum distance between $t_j$-dimensional subspaces of $\bbF_q^n$. Moreover, $\cF_j$ and $\cF'_j$ are different and we obtain $|\cC_i|\leq |\cC_j|$, for all $i\leq j\leq r$.
\end{enumerate}
\end{proof}

\begin{remark}
Notice that having a projected code attaining the maximum distance it is not enough to deduce that other projected codes satisfy the same property. The following example shows this situation. 
\end{remark}

\begin{example}
Consider the canonical basis $\{\be_1,\dots, \be_6\}$ of $\bbF_q^6$ and let $\cC$ be the flag code of type $(2, 3)$ on $\bbF_q^6$ given by the flags
$$
\begin{array}{cccccc}
\cF^1  & = & ( \langle \be_1, \be_2 \rangle, & \langle \be_1, \be_2, \be_3 \rangle), \\
\cF^2  & = & ( \langle \be_1, \be_3 \rangle, & \langle \be_1, \be_2, \be_3 \rangle),\\
\cF^3  & = & ( \langle \be_4, \be_5 \rangle, & \langle \be_4, \be_5, \be_6 \rangle). \\
\end{array}
$$
Notice that $\cC_2= \{\langle \be_1, \be_2, \be_3 \rangle, \langle \be_4, \be_5, \be_6 \rangle\}$ is a code of maximum distance, $d_S(\cC_2)=6$. However, the first projected code $\cC_1=\{ \langle \be_1, \be_2 \rangle, \langle \be_1, \be_3 \rangle,  \langle \be_4, \be_5 \rangle \}$ has distance $d_S(\cC_1)=  d_S( \cF_1^1,  \cF_1^2) = 2$, whereas the maximum distance for its dimension is $4$. 
\end{example}

We will use Proposition \ref{prop: a projected code with max dist} in order to characterize ODFCs in terms of, at most, two of their projected codes. These constant dimension codes are the ones (if they exist) of closest dimensions  to $\frac{n}{2}$  in the type vector, both at left and right. Let us make this idea precise. Given an arbitrary but fixed type vector $(t_1, \dots, t_r)$ and an ambient space $\bbF_q^n$, we give special attention to two indices  defined as
\begin{equation}\label{eq: indices a and b}
\begin{array}{ccll}
a &=& \max \{ i\in \{1, \dots, r\} \ | \ 2t_i\leq n \}   & \text{and}\\
b &=& \min \{ i\in \{1, \dots, r\} \ | \ 2t_i \geq n\}. &
\end{array}
\end{equation}
Note that the sets 
\[
\{ i\in \{1, \dots, r\} \ | \ 2t_i\leq n\} \ \text{and} \ \{ i\in \{1, \dots, r\} \ | \ 2t_i \geq n\}
\] 
cover the family of indices $\{1, \dots, r\}$. Hence, at least one of them must be nonempty. Even more, the first set is empty if, and only if, $b=1$ and every dimension in the type vector is lower bounded by $\frac{n}{2}$. Similarly, the second one does not contain any element if, and only if, $a=r$, i.e., if all the dimensions in the type vector are upper bounded by $\frac{n}{2}$. In any other situation, $a$ and $b$ are well-defined and $a\leq b$. The equality holds if, and only if, $n$ is even and the dimension $\frac{n}{2}$ appears in the type vector. In this case $t_a=t_b=\frac{n}{2}$. In any other situation, these two sets partition 
\[
\{1, \dots, r\}=\{ 1, \dots, a\} \ \dot\cup \ \{ b, \dots,  r \},
\] 
with $b=a+1$.

For sake of simplicity, the following results are presented in terms of both indices $a$ and $b$. Despite the fact that these two indices do not need to exist simultaneously, at least one of them is always well defined. In that way, for type vectors in which the index $a$ (resp. $b$) does not exist, the next result still holds true and gives a characterization of ODFCs just in terms of the projected code of dimension $t_b=t_1$ (resp. $t_a=t_r$). In any case, it represents a significant improvement with respect to Theorem \ref{teo:carac_odfc} (see \cite[Th. 3.11]{CasoPlanar}).
 
\begin{theorem}\label{teo: characterization ODFC a and b}
Let $\cC$ be a flag code of type $(t_1, \dots, t_r)$ on $\bbF_q^n$ and consider indices $a$ and $b$ as in (\ref{eq: indices a and b}). The following statements are equivalent:
\begin{enumerate}[$(i)$]
\item The flag code $\cC$ is an ODFC.
\item $\cC_{a}$ and $\cC_{b}$ are constant dimensions codes of maximum distance with cardinality $|\cC_{a}|=|\cC_{b}|=|\cC|$.
\end{enumerate} 
\end{theorem}
\begin{proof} 
The implication $(i) \then (ii)$ follows straightforwardly from Theorem \ref{teo:carac_odfc}. To show 
$(ii) \then (i)$, notice that if $\cC_{a}$ and $\cC_{b}$ have the maximum possible distance, by means of Proposition \ref{prop: a projected code with max dist}, we have
\begin{eqnarray*}
 |\cC| & = & |\cC_{a}|\leq |\cC_i| \  \text{for every} \ i\leq a \ \text{and} \\
 |\cC| & = & |\cC_{b}|\leq |\cC_j| \   \text{for every} \ j\geq b. 
\end{eqnarray*}
Since the cardinality of every projected code is upper bounded by $|\cC|$, we conclude that $|\cC|=|\cC_1|=\dots=|\cC_r|$, i.e., the flag code $\cC$ is disjoint. Now, in order to see that every projected code attains the maximum possible distance, we argue as follows. Take an index $1\leq i\leq r$ and consider a pair of different subspaces $\cF_i, \cF'_i\in\cC_i$. These subspaces come from different flags $\cF, \cF'\in\cC$  and, since $\cC$ is disjoint, we have $\cF_j\neq\cF'_j$ for every $1\leq j\leq r$. In particular, $\cF_{a}\neq \cF'_{a}$ and $\cF_{b}\neq \cF'_{b}$. Hence, by means of Proposition \ref{prop: a projected code with max dist}, in any case, the distance $d_S(\cF_i, \cF'_i)$ is the maximum possible one for dimension $t_i$ and, as a result, every projected code $\cC_i$ is a constant dimension code of maximum distance. Thus, by application of Theorem \ref{teo:carac_odfc}, the flag code $\cC$ is an ODFC.
\end{proof}

Here below we translate the previous results into the orbital scenario. As before, we express them in terms of both indices $a$ and $b$ defined in (\ref{eq: indices a and b}), but always having in mind that one of them might not exist. In such a case, the result is fulfilled by the remaining index.

The following proposition was already stated in \cite{PlanarOrbConstr} for $n=2k$ and flags of full type vector. 

\begin{proposition}\label{prop: stab_inclusion}
Let $\cF$ be a flag of type $(t_1, \dots, t_r)$ on $\bbF_q^n$ and $\bH$ a subgroup of $\GL(n,q)$ such that $\orb_{\bH}(\cF_a)$ and $\orb_{\bH}(\cF_{b})$ are subspace codes of maximum distance. Then
\begin{enumerate}[$(i)$]
\item $\stab_{\bH}(\cF_i)\subseteq \stab_{\bH}(\cF_a)$ for all $i\leq a$ and 
\item $\stab_{\bH}(\cF_i)\subseteq \stab_{\bH}(\cF_{b})$ for all $i\geq b$.
\end{enumerate}
\end{proposition}
\begin{proof}
 $(i)$ Consider a matrix $A\in \bH\setminus\stab_{\bH}(\cF_a)$. Then $\cF_a\neq\cF_a\cdot A$ and, by means of Proposition \ref{prop: a projected code with max dist}, for every $i\leq a$, it holds $d_S(\cF_i, \cF_i \cdot A)=2t_i$. In particular, it is clear that $A\not\in \stab_{\bH}(\cF_i),$ whenever $i\leq a$. Equivalently, we have that $\stab_{\bH}(\cF_i)\subseteq \stab_{\bH}(\cF_a)$ for all $i\leq a$.

\noindent $(ii)$ If $A\in\bH\setminus\stab_{\bH}(\cF_{b})$, then $\cF_{b}\neq\cF_{b}\cdot A$ and Proposition \ref{prop: a projected code with max dist} leads to $d_S(\cF_i, \cF_i \cdot A)=2(n-t_i)$ for all $i\geq b.$ Consequently, $A\not\in \stab_{\bH}(\cF_i)$, for every $i\geq b$, and the result holds. 
\end{proof}

\begin{remark}
Observe that Proposition \ref{prop: a projected code with max dist} gives us some conditions on the cardinality of the projected codes. From that result, in the orbital scenario, if we assume that $\cC_a$ and $\cC_b$ attain the maximum possible distance, it follows that the order of every $\stab_{\bH}(\cF_i)$ must be upper bounded either by $|\stab_{\bH}(\cF_a)|$ or  $|\stab_{\bH}(\cF_{b})|$. However, in Proposition \ref{prop: stab_inclusion}, we obtain a stronger condition: a subgroup relationship.
\end{remark}

Next we summarize several different characterizations for ODFCs arising as orbits of the action of an arbitrary subgroup of $\GL(n,q)$.

\begin{theorem}\label{teo:cond_suf_ODFC_orbital}
Consider a flag $\cF$ of type $(t_1, \dots, t_r)$ on $\bbF_q^n$ and a subgroup $\bH$ of $\GL(n,q).$ The following statements are equivalent:
\begin{enumerate}[$(i)$]
\item The flag code $\orb_{\bH}(\cF)$ is an ODFC. 
\item The subspace codes $\orb_{\bH}(\cF_a)$ and $\orb_{\bH}(\cF_{b})$ are of maximum distance and $\stab_{\bH}(\cF_a)=\stab_{\bH}(\cF_{b})\subseteq\stab_{\bH}(\cF_i),$ for every $1\leq i\leq r$. 
\item The subspace codes $\orb_{\bH}(\cF_a)$ and $\orb_{\bH}(\cF_{b})$ are of maximum distance and $\stab_{\bH}(\cF_a)=\stab_{\bH}(\cF_{b})\subseteq\stab_{\bH}(\cF)$.
\item The subspace codes $\orb_{\bH}(\cF_a)$ and $\orb_{\bH}(\cF_{b})$ are of maximum distance and $|\stab_{\bH}(\cF_a)|=|\stab_{\bH}(\cF_{b})|\leq|\stab_{\bH}(\cF)|$.
\end{enumerate}
The cardinality of such a flag code is $|\orb_{\bH}(\cF)|=|\orb_{\bH}(\cF_a)|=|\orb_{\bH}(\cF_b)|.$
\end{theorem}
\begin{proof}
Observe that, by means of Theorem \ref{teo:carac_odfc}, together with Proposition \ref{prop:disjorb}, statement $(i)$ clearly implies the other ones. On the other hand, by application of Proposition \ref{prop: stab_inclusion} and expression (\ref{eq:inter_stab}), all conditions $(ii)$, $(iii)$ and $(iv)$ are equivalent and make the code $\orb_{\bH}(\cF)$ be disjoint by Proposition \ref{prop:disjorb}. In particular, it holds $|\orb_{\bH}(\cF)|=|\orb_{\bH}(\cF_a)|=|\orb_{\bH}(\cF_b)|$ with projected codes of dimensions $t_a$ and $t_b$ attaining the maximum possible distances. Theorem \ref{teo: characterization ODFC a and b} finishes the proof. 
\end{proof}

Recall that the cardinality of an orbit flag code is completely determined by the orders of the acting group and the stabilizer subgroup of the generating flag. Fixed the acting group, a natural way of obtaining codes with better cardinalities is to consider the union of different orbits.  We finish this section by characterizing when the union  of orbit flag codes is an ODFC. To this purpose, notice that every nonempty subset of an ODFC is either an ODFC too or a trivial code having just one element. In both cases, such a subset is a disjoint flag code. With the goal of obtaining better cardinalities in mind, we proof first the following lemma where we work with the union of two disjoint flag codes that are orbits of the same group.

\begin{lemma}\label{lemma: union 2 orbits}
Let $\cF$ and $\cF'$ be flags of type $(t_1, \dots, t_r)$ on $\bbF_q^n$ and take a subgroup $\bH$ of $\GL(n, q)$ such that the orbits $\orb_\bH(\cF)$ and $\orb_\bH(\cF')$ are disjoint flag codes. The following statements hold:
\begin{enumerate}[$(a)$]
\item  If $\orb_\bH(\cF_i)\neq \orb_\bH(\cF'_i)$, for some $i\in\{1,\dots, r\}$, then $|\orb_\bH(\cF)\cup \orb_\bH(\cF')|=|\orb_\bH(\cF)|+|\orb_\bH(\cF')|$.
\item If $\orb_\bH(\cF_i)\neq \orb_\bH(\cF'_i)$, for all $i=1,\dots, r$, then $\orb_\bH(\cF)\cup\orb_\bH(\cF')$ is a disjoint flag code.
\item If $\orb_\bH(\cF_i)=\orb_\bH(\cF'_i)$, for some $i\in\{1,\dots, r\}$, then $\orb_\bH(\cF)\cup\orb_\bH(\cF')$ is a disjoint flag code if, and only if, $\orb_\bH(\cF)=\orb_\bH(\cF')$.
\end{enumerate} 
\end{lemma}
\begin{proof}
$(a)$ Since  $\orb_\bH(\cF_i)\neq \orb_\bH(\cF'_i)$,  their intersection is the empty set and, since $\orb_\bH(\cF)$ and $\orb_\bH(\cF')$ are disjoint flag codes, it follows that
\begin{equation}\label{eq:orb_disj}
\begin{array}{ccl}
|\orb_\bH(\cF)\cup\orb_\bH(\cF')| & \leq & |\orb_\bH(\cF)|+ |\orb_\bH(\cF')| =  |\orb_\bH(\cF_i)|+ |\orb_\bH(\cF'_i)|\\
                                                     & =     & |\orb_\bH(\cF_i)\cup \orb_\bH(\cF'_i)|  =  |(\orb_\bH(\cF)\cup \orb_\bH(\cF'))_i|\\
                                                     & \leq & |\orb_\bH(\cF)\cup\orb_\bH(\cF')|.
\end{array}
\end{equation}
Therefore, $|\orb_\bH(\cF)\cup\orb_\bH(\cF')| = |\orb_\bH(\cF)|+ |\orb_\bH(\cF')| $ and the statement holds.

$(b)$ Follows directly from applying $(\ref{eq:orb_disj})$ for all $i=1,\dots, r$.

$(c)$
Assume that the union code $\orb_\bH(\cF)\cup\orb_\bH(\cF')$ is a disjoint flag code. Given that $\orb_\bH(\cF_i)=\orb_\bH(\cF'_i)$, there must exist a matrix $A\in\bH$ such that $\cF_i=\cF'_i\cdot A$. Hence, both flags $\cF$ and $\cF'\cdot A$ are in $\orb_\bH(\cF)\cup\orb_\bH(\cF')$ and they share their $i$-th subspace. Thus, we conclude that $\cF=\cF'\cdot A$ and therefore $\orb_\bH(\cF)=\orb_\bH(\cF')$. The converse statement trivially holds.
\end{proof}

Now, with the benefit of the above lemma and using Theorem \ref{teo: characterization ODFC a and b}, the next result states when the union of a family of disjoint flag codes, arising as orbits of the same group, provides ODFCs of larger cardinality. To do so, we consider indices $a$ and $b$ as in (\ref{eq: indices a and b}).

\begin{theorem}\label{theo:unionODFCorbitals}
Let $\{\cF^j=(\cF^j_1,\dots, \cF^j_{r})\}_{j=1}^m$ be a family of flags of type $(t_1, \dots, t_r)$ on $\bbF_q^n$ and consider a subgroup $\bH$ of $\GL(n,q)$ such that every orbit $\orb_{\bH}(\cF^j)$ is a disjoint flag code, for $1\leq j\leq m$. If the subspaces $\cF^1_{a}, \dots  ,\cF^m_{a}$ and $\cF^1_{b}, \dots  ,\cF^m_{b}$ lie in different orbits under the action of $\bH$, then
\[
|\bigcup_{j=1}^m \orb_{\bH}(\cF^j)| =  \sum_{j=1}^m |\orb_{\bH}(\cF^j)|.
\]

Moreover, the following statements are equivalent:
\begin{enumerate}[$(i)$]
\item The union flag code $\bigcup_{j=1}^m \orb_{\bH}(\cF^j)$ is an ODFC.
\item The projected union codes $\bigcup_{j=1}^m \orb_{\bH}(\cF^j_a)$ and $\bigcup_{j=1}^m \orb_{\bH}(\cF^j_{b})$ have the maximum possible distance.
\end{enumerate} 
\end{theorem}
\begin{proof}
Since every orbit $\orb_\bH(\cF^j)$ is  a disjoint flag code, we can argue as in $(\ref{eq:orb_disj})$ of the previous lemma to obtain that 
\begin{equation}\label{eq:sum_orbdisj}
|\bigcup_{j=1}^m \orb_{\bH}(\cF^j)| = \sum_{j=1}^m |\orb_{\bH}(\cF^j)|  = \sum_{j=1}^m |\orb_{\bH}(\cF^j_i)|  ,                                      
\end{equation}
for $i=a, b$. Hence, by means of Theorem \ref{teo: characterization ODFC a and b}, the union flag code is an ODFC if, and only if, its projected codes of dimensions $t_a$ and $t_b$ attain the maximum possible distance.
\end{proof}

We will use these theoretical results in the following section in order to give specific constructions of orbit ODFCs having the maximum possible cardinalities for the corresponding type vectors.


\section{ODFC From Singer Groups}\label{sec:construc}

This section is devoted to construct flag codes of maximum distance having an orbital structure. For this, we will use suitable Singer groups (or their subgroups) and their transitive action on lines and hyperplanes (Theorem \ref{teo:stabSinger_rectes}). The goal of Section \ref{subsect:$n=ks$} is to obtain ODFCs on $\bbF_q^n$ having a $k$-spread as a projected code, for $k$ a divisor of $n$. To do so, according to Theorem \ref{teo:max_ad_type}, we consider first flags of type $(1,\dots, k,n-k,\ldots, n-1)$. Such a construction leads to full flag codes whenever $n=2k$ or $k=1$ and $n=3$. In Section \ref{subsect:$n=2t+1$}, we build ODFCs of full type vector for the remaining cases.


\subsection{Orbit ODFC from Desarguesian spreads}\label{subsect:$n=ks$}

In this section we address the orbital construction of ODFCs on $\bbF_q^{n}$ with the $k$-spread $\cS$ defined in (\ref{def:spreadS}) as a projected code. To this end, write $n=ks$ for some $s\geq 2$. Recall that, by virtue of Theorem \ref{teo:ODFC_maxcard}, such a code has also the largest possible size. 
Throughout the rest of this section, and for sake of simplicity, we will write $\cF=(\cF_1,\dots,\cF_k,\cF_{n-k},\dots,\cF_{n-1})$ to denote an arbitrary flag of the {\em full admissible type vector} $(1, \dots, k, n-k, \dots, n-1)$.  
In these conditions, indices $a$ and $b$ of (\ref{eq: indices a and b}) are such that $t_a=k$ and $t_b=n-k$. Consider the Singer group $\langle \psi(M_s)\rangle$ of $\GL(n,q)$ defined in Section \ref{sec:spread_Singer}. Recall that the constant dimension codes $\cS$ and $\cH$ arise as their orbits (see (\ref{eq:S_orbSinger}) and (\ref{eq:H_orbSinger})).  Next, we use Theorem \ref{teo:cond_suf_ODFC_orbital} in order to characterize those subgroups of of $\langle \psi(M_s)\rangle$ that are appropriate to construct ODFCs.

\begin{theorem}\label{teo:orb_ODFC_singer}
Let $\cF=(\cF_1,\dots,\cF_k,\cF_{n-k},\dots,\cF_{n-1})$ be a flag of full admissible type vector such that $\cF_{k}\in\cS$ and $\cF_{n-k}\in \cH$. For any positive integer $t$ dividing $q^{n}-1$, consider the unique subgroup $\bT$ of $\langle \psi(M_s)\rangle$ of order $t$. Then:
\begin{enumerate}[$(i)$]
\item $|\orb_{\bT}(\cF)|=\frac{t}{\gcd(t,q-1)}$.
\item $\orb_{\bT}(\cF)$ is  an ODFC if, and only if,  $\gcd(t,q^k-1)=\gcd(t,q-1)\neq t$.
\end{enumerate}
\end{theorem}
\begin{proof}
 $(i)$ By $(\ref{eq:inter_stab})$ and  Theorem \ref{teo:stabSinger_rectes}, it follows that 
\[
\stab_{\langle \psi(M_s)\rangle}(\cF)=\stab_{\langle \psi(M_s)\rangle}(\cF_1)=\stab_{\langle \psi(M_s)\rangle}(\cF_{n-1})=\{aI_n\ |\ a\in\bbF_q^*\}, 
\]
which has order $q-1$. As a result, $\stab_{\bT}(\cF)  =  \bT\cap \{aI_n\ |\ a\in\bbF_q^*\}$, is a group of order 
$\gcd(t,q-1)$ and the statement holds. Notice that, in particular, $\orb_\bT(\cF)=\{\cF\}$ exactly when $t\mid q-1$.

$(ii)$ By using (\ref{eq:stab_G}),  it follows that 
\begin{eqnarray}
\stab_{\bT}(\cF_k) & = & \bT\cap \stab_{\langle \psi(M_s)\rangle}(\cF_k)=\bT\cap \{\psi(aI_s)\ |\ a\in\bbF_{q^k}^*\}= \nonumber \\
                              & = &  \bT\cap \stab_{\langle \psi(M_s)\rangle}(\cF_{n-k})=\stab_{\bT}(\cF_{n-k}), \nonumber
\end{eqnarray}
which is a group of order $\gcd(t, q^k-1)$. Since we are working with cyclic groups, we obtain that $\stab_\bT(\cF_k)=\stab_{\bT}(\cF_{n-k})= \stab_\bT(\cF)$ if, and only if, $\gcd(t, q^k-1)=\gcd(t, q-1)$. On the other hand, since $\orb_\bT(\cF_k)\subseteq \cS$ and $\orb_\bT(\cF_{n-k})\subseteq\cH$, these projected codes of the flag code $\orb_{\bT}(\cF)$  will be constant dimension codes of maximum distance whenever they have at least $2$ elements. As a result, the statement $(ii)$ follows from  Theorem \ref{teo:cond_suf_ODFC_orbital} and $(i)$.
\end{proof}

As stated in Theorem \ref{teo:ODFC_maxcard}, the size of the $k$-spread $\cS$, that is, $\frac{q^n-1}{q^k-1}$, determines the maximum size of an ODFC having $\cS$ as a projected code. In order to achieve this optimal size, we will consider unions of orbits under the action of a suitable subgroup $\bT$ of $\langle \psi(M_s)\rangle$ and then apply Theorem \ref{theo:unionODFCorbitals} as follows. 

\begin{theorem}\label{teo:unio_orb_ODFC_singer}
Let $\bT$ be a subgroup of order $t$ of $\langle \psi(M_s)\rangle$ such that $\gcd(t, q^k-1)=\gcd(t,q-1)$. For $m\geq 2$, let $\{\cF^j=(\cF^j_1,\dots, \cF^j_{k},\cF^j_{n-k},\dots, \cF^j_{n-1})\}_{j=1}^m$ be a family of flags of full admissible type such that $\cF^1_{k},\dots,\cF^m_{k}\in\cS$ and  $\cF^1_{n-k},\dots,\cF^m_{n-k}\in\cH$ lie in different orbits under the action of $\bT$. Then
\begin{enumerate}[$(i)$]
\item $\cup_{j=1}^m \orb_{\bT}(\cF^j)$ is an ODFC of size $\frac{m t}{\gcd(t, q-1)}$. 
\item If $m= \frac{(q^n-1)\gcd(t, q-1)}{(q^k-1) t}$, then $\cup_{j=1}^m \orb_{\bT}(\cF^j)$ is an ODFC of the maximum size, that is, $\frac{q^n-1}{q^k-1}$.
\end{enumerate}
\end{theorem}
\begin{proof}
$(i)$ By Theorem \ref{teo:orb_ODFC_singer} we know that $|\orb_{\bT}(\cF^j)|=\frac{t}{\gcd(t,q-1)}$, for all $1\leq j\leq m$. Moreover, since $\gcd(t, q^k-1)=\gcd(t,q-1)$, the same theorem ensures that either $\orb_{\bT}(\cF^j)=\{\cF^j\}$ or $\orb_{\bT}(\cF^j)$ is an ODFC, for all $1\leq j\leq m$. In any case, one has that $\orb_{\bT}(\cF^j)$ is a disjoint flag code, for all $1\leq j\leq m$. Thus, we have the hypotheses of Theorem \ref{theo:unionODFCorbitals} and following $(\ref{eq:sum_orbdisj})$ we can write
\[
|\cup_{j=1}^m \orb_{\bT}(\cF^j)|= \sum_{j=1}^m |\orb_{\bT}(\cF^j)| = \sum_{j=1}^m |\orb_{\bT}(\cF^j_i)| =  \frac{m t}{\gcd(t, q-1)} \geq 2,
\]
for $i=k, n-k$. In particular, the projected codes of dimensions $k$ and $n-k$ of the union flag code $\cup_{j=1}^m \orb_{\bT}(\cF^j)$ are subsets of $\cS$ and $\cH$, respectively, having at least two elements. Hence they are subspace codes of maximum distance and Theorem \ref{theo:unionODFCorbitals} states that $\cup_{j=1}^m \orb_{\bT}(\cF^j)$  is an ODFC.

Statement $(ii)$ follows just by computing the number of orbits of the action of $\bT$ on $\cS$. 
\end{proof}

\begin{remark}
Theorem \ref{teo:orb_ODFC_singer} states which subgroups of  $\langle \psi(M_s)\rangle$  allow the construction of ODFCs as a single orbit of them. Notice that bigger subgroups not always will provide bigger orbit flag codes. In addition, it may eventually happen that some subgroup provides an orbit ODFC of the maximum possible size,  $\frac{q^n-1}{q^k-1}$. Otherwise, Theorem \ref{teo:unio_orb_ODFC_singer} leads to an optimal construction consisting of the union of several orbits. Clearly, the larger the size of each orbit, the fewer orbits we need to join to reach the maximum size and vice versa. In particular, the degenerate case where an ODFC is constructed as a union of $\frac{q^n-1}{q^k-1}$ orbits with just one element is also contemplated in Theorem  \ref{teo:unio_orb_ODFC_singer}. 

All these considerations are reflected in the following examples, in which we apply Theorems \ref{teo:orb_ODFC_singer} and \ref{teo:unio_orb_ODFC_singer} for different values of the parameters. 
\end{remark}
 
\begin{examples} 
With the notation of Theorem \ref{teo:orb_ODFC_singer}, we consider all the divisors $t$ of $q^n-1$ such that $\gcd(t, q^k-1)=\gcd(t,q-1)$ and  the corresponding subgroup $\bT$ of $\langle\psi(M_k)\rangle$ of order $t$. Consider a flag $\cF$ of the full admissible type $(1,\dots, k,n-k,\dots, n-1)$ such that $\cF_k\in\cS$ and $\cF_{n-k}\in\cH$. Finally, denote by $m$ the number of required orbits of $\bT$ to attain the maximum size,  $\frac{q^n-1}{q^k-1}$, for an ODFC with these parameters. 
\begin{enumerate}
\item Put $q=3$, $k=3$ and $n=6$. Thus, $k=n-k$, $q^n-1=728$, $q^k-1=26$ and $\frac{q^n-1}{q^k-1}=28$. Then  
\begin{table}[H]
\centering
\begin{tabular}{|c|c|c|c|c|c|c|c|c|}
\hline
$t$                         & 1 & 2 & 4 & 7 & 8 & 14  & 28 & 56  \\ \hline
$|\orb_\bT(\cF)|$           & 1 & 1 & 2 & 7 & 4 &  7  & 14 & 28 \\ \hline  
$m$                         & 28 & 28  & 14 & 4 & 7  & 4  & 2 & 1    \\ \hline
\end{tabular}
\caption{$q=3$, $k=3$ and $n=6$.}
\end{table}
Notice that, in this case, the subgroup of order $t=56$ allows us to obtain ODFCs of full type vector and having the best possible size, i.e., $28$, by using a single orbit. In this sense, for odd characteristic, Theorem \ref{teo:orb_ODFC_singer} eventually improves the construction presented in \cite[Prop. 4.15]{PlanarOrbConstr}, where two orbits were always needed. Moreover, remark that the subgroup of order $t=8$ gives an orbit ODFC of smaller size than the obtained with the subgroup of order $t=7$.

\item Put $q=4$, $k=3$ and $n=9$. Thus, $n-k=6$, $q^n-1=262143$, $q^k-1=63$ and $\frac{q^n-1}{q^k-1}=4161$. Then
\begin{table}[H]
\centering
\begin{tabular}{|c|c|c|c|c|c|c|c|c|}
\hline
$t$                         & 1   & 3   & 19 & 57 & 73 & 219 & 1387 & 4161  \\ \hline
$|\orb_\bT(\cF)|$           & 1   & 1   & 19 & 19 & 73 &  73 & 1387 & 1387          \\ \hline  
$m$                         & 4161 & 4161 & 219 & 219 & 57 & 57 & 3 &   3 \\ \hline 
\end{tabular}
\caption{$q=4$, $k=3$ and $n=9$.}
\end{table}
The largest orbit size is $1387$ and it is obtained when the acting group has order either $1387$ or $4161$. On the other hand, the maximum possible size of an ODFC with these parameters is $4161$. Hence, in order to achieve that cardinality, we must consider the union of, at least, $3$ different orbits.
\end{enumerate}
\end{examples}

The orbital constructions of ODFC provided in this section present a restriction on the type vector, coming from the condition of having a spread as a projected code. However, there are two possible situations in which flag codes of full type can be given by using Theorems \ref{teo:orb_ODFC_singer} and \ref{teo:unio_orb_ODFC_singer}. First, for even values of $n$, taking the divisor $k=\frac{n}{2}$ leads to a construction of full type in which the values $k$ and $n-k$ coincide. This particular case was first studied in \cite{PlanarOrbConstr}, where a construction using the action of a Singer subgroup of $\SL(2k, q)$ is presented. On the other hand, for $n=3$ and $k=1$, the action of a Singer subgroup of $\GL(3,q)$ on the Grassmannian of lines and hyperplanes gives us a construction of type $(1,2)$, this construction is also known and the reader can find it in \cite[Prop. 2.5]{Kurz20}, where the author shows that it is the one with the biggest cardinality among ODFC of full type when $n=3$. In the following section, we consider the remaining situations, that is, we address the construction of orbit ODFCs of full type vector on $\bbF_q^n$ for odd values of $n> 3$.


\subsection{Orbit ODFC of full type vector}\label{subsect:$n=2t+1$}

Throughout this section, we work with full flags on $\bbF_q^{2k+1}$, for some $k>1$. In this case, by virtue of Theorem \ref{teo: characterization ODFC a and b}, the construction of ODFCs can be done by giving appropriate constant dimension codes for dimensions $k$ and $k+1$ (see (\ref{eq: indices a and b})). To this end, we present the next subgroup of $\GL(2k+1, q)$. Let $M_{k+1}\in \GL(k+1, q)$ be the companion matrix of a primitive polynomial of degree $k+1$ in $\bbF_q[x]$. Recall that, as pointed out in Section \ref{subsec:subspaces}, $M_{k+1}$ is a Singer cycle of $\GL(k+1, q)$ and $\bbF_q[M_{k+1}]=\langle M_{k+1}\rangle\cup\{0_{(k+1)\times (k+1)}\}$ is a matrix representation of the finite field of $q^{k+1}$ elements.  Let us write 
\[
g=\begin{pmatrix}
\begin{array}{c|c}
I_k & 0_{k\times(k+1)}\\ \hline
0_{(k+1)\times k} & M_{k+1}
\end{array}
\end{pmatrix}
\in \GL(2k+1, q)
\]
and consider the cyclic group 
\begin{equation}\label{eq: group 2t+1}
\bG= \left\langle 
g
\right\rangle
= \left \lbrace
g^i \ | \ 0\leq i\leq q^{k+1}-2
\right\rbrace .
\end{equation}
Clearly, $\bG$ is a subgroup of order $q^{k+1}-1$ of $\GL(2k+1, q)$, isomorphic to the Singer subgroup $\langle M_{k+1}\rangle$ of $\GL(k+1, q)$. In the rest of this section, the orbit codes considered will be always generated by the action of this particular group $\bG$.

We start by characterizing the subspaces of dimensions $k$ and $k+1$ of $\bbF_q^{2k+1}$ whose orbits under the action of $\bG$ are constant dimension codes of maximum distance. Given arbitrary subspaces $\cU=\rsp(U)\in\cG_q(k, 2k+1)$ and $\cV=\rsp(V)\in \cG_q(k+1, 2k+1)$, the  respective full-rank generator matrices $U\in\bbF_q^{k\times(2k+1)}$ and $V\in\bbF_q^{(k+1)\times(2k+1)}$ can split into two blocks as 
\begin{equation}\label{eq: n=2k+1 generatos matrices}
U=( U_1 \ | \ U_2 ) \ \ \text{and} \ \ V=( V_1 \ | \ V_2 ),
\end{equation}
where $U_1$ (resp. $V_1$) denotes the first $k$ columns of $U$ (resp. $V$). Therefore, $U_1\in\bbF_q^{k\times k}$, $U_2\in\bbF_q^{k\times(k+1)}$, $V_1\in\bbF_q^{(k+1)\times k}$ and $V_2\in\bbF_q^{(k+1)\times(k+1)}$. Using this notation, we can write
\begin{eqnarray}\label{eq: subspaceU in projected codes n=2k+1}
\orb_\bG(\cU) & = & \{\cU \cdot g^i\ |\ 0\leq i\leq q^{k+1}-2\} \nonumber \\
              & = & \{\rsp( U_1 \ | \ U_2 M_{k+1}^i )  \ |\ 0\leq i\leq q^{k+1}-2\}
\end{eqnarray}
and
\begin{eqnarray}\label{eq: subspaceV in projected codes n=2k+1}
\orb_\bG(\cV) & = & \{\cV \cdot g^i \ |\ 0\leq i\leq q^{k+1}-2\} \nonumber \\
              & = & \{\rsp( V_1 \ | \ V_2 M_{k+1}^i ) \ |\ 0\leq i\leq q^{k+1}-2\}.
\end{eqnarray}

With this notation, the following results hold. To make this section easier to read, their proofs, and all the ones concerning subspace codes, are included in the final Appendix of the article, so that only proofs concerning results on flag codes appear here. 

\begin{proposition}\label{prop: n=2k+1 projected code k}
The orbit code $\orb_\bG(\cU)$ defined in (\ref{eq: subspaceU in projected codes n=2k+1}) is a partial spread of dimension $k$ of $\bbF_q^{2k+1}$ if, and only if, $\rk(U_1)=\rk(U_2)=k$. Its cardinality is  $|\orb_\bG(\cU)|=|\bG|=q^{k+1}-1.$
\end{proposition}

\begin{proposition}\label{prop: n=2k+1 projected code k+1}
The orbit code $\orb_\bG(\cV)$ defined in (\ref{eq: subspaceV in projected codes n=2k+1}) attains the maximum possible distance if, and only if, $\rk(V_1)=k$ and $\rk(V_2)=k+1$. Its size is  $|\orb_\bG(\cV)|=|\bG|=q^{k+1}-1.$
\end{proposition}

Here below, we use the previous characterizations for constant dimension codes of maximum distance in order to provide orbit ODFCs of full type on $\bbF_q^{2k+1}$. To do so, we need to consider nested subspaces $\cU\subsetneq \cV$ of dimensions $k$ and $k+1$, respectively. Using the notation of (\ref{eq: n=2k+1 generatos matrices}),  we can formulate the problem in a matrix approach: given a full-rank generator matrix $U=(U_1 \ |\  U_2)\in\bbF_q^{k\times(2k+1)}$ of $\cU$, we consider a subspace $\cV$ spanned by the rows of a matrix  $V\in\bbF_q^{(k+1)\times(2k+1)}$, obtained by adding an appropriate row to $U$. In other words, we choose vectors $\bv_1\in\bbF_q^{k}$ and $\bv_2\in\bbF_q^{k+1}$ such that the matrix

\begin{equation}\label{eq: generator matrix V n=2k+1}
V= ( V_1 \ | \ V_2) =
\begin{pmatrix}
\begin{array}{c|c}
U_1   & U_2\\ \hline
\bv_{1} &  \bv_{2}\\
\end{array}
\end{pmatrix}
\end{equation}
has rank equal to $k+1$. Using this notation, we present the next construction of ODFCs arising from the action of the group $\bG$ defined in (\ref{eq: group 2t+1}).

\begin{theorem}\label{theo: n=2k+1 ODFC orbita de G}
Let $\cF=(\cF_1,\dots,\cF_{2k})$ be a full flag on $\bbF_q^{2k+1}$ such that 
\[
\cF_k=\cU=\rsp(U_1 \ | \ U_2) \ \text{and} \ \cF_{k+1}=\cV=\rsp(V_1 \ |  \ V_2),
\]
with generator matrix $(V_1 \ | \ V_2)$ as in (\ref{eq: generator matrix V n=2k+1}) and consider the group $\bG$ defined in (\ref{eq: group 2t+1}). The following statements are equivalent:
\begin{enumerate}[(i)]
\item the flag code $\orb_\bG(\cF)$ is an ODFC.
\item  $U_1\in\bbF_q^{k\times k}$ and $V_2\in\bbF_q^{(k+1)\times(k+1)}$ are invertible matrices.
\end{enumerate}
In this situation, $|\orb_\bG(\cF)|=|\bG|=q^{k+1}-1$.
\end{theorem}
\begin{proof}
Assume that $\orb_\bG(\cF)$ is an ODFC. In particular, the projected codes $\orb(\cF_{k})$ and $\orb_\bG(\cF_{k+1})$ must attain the maximum distance. By means of Propositions \ref{prop: n=2k+1 projected code k} and \ref{prop: n=2k+1 projected code k+1}, it must hold $\rk(U_1)=\rk(U_2)=\rk(V_1)=k$ and $\rk(V_2)=k+1$. Consequently, $U_1$ and $V_2$ are invertible matrices.

Conversely, assume now that $\rk(U_1)=k$ and $\rk(V_2)=k+1$. Since $U_2$ is composed by the first $k$ rows of the invertible matrix $V_2$, we clearly obtain that $\rk(U_2)=k$. On the other hand, observe that $V_1\in\bbF_q^{(k+1)\times k}$ contains $U_1$ as a submatrix. Hence, its rank is $k$ as well. Now, by using Propositions \ref{prop: n=2k+1 projected code k} and \ref{prop: n=2k+1 projected code k+1}, we conclude that both $\orb_\bG(\cF_k)$ and $\orb_\bG(\cF_{k+1})$ are constant dimension codes of maximum distance such that $\stab_\bG(\cF_k)=\stab_\bG(\cF_{k+1})=\{ I_{2k+1}\}\subseteq \stab_\bG(\cF_i)$, for every $1\leq i\leq 2k$. Hence, by application of Theorem \ref{teo:cond_suf_ODFC_orbital}, the orbit flag code $\orb_\bG(\cF)$ is an ODFC of size $|\orb_\bG(\cF)|=|\bG|=q^{k+1}-1$. 
\end{proof}

The ODFC constructed in Theorem \ref{theo: n=2k+1 ODFC orbita de G} contains $q^{k+1}-1$ flags. It is the largest size for orbits under the action of the group $\bG$.  On the other hand, as proved in \cite[Prop. 2.4]{Kurz20}, the maximum possible cardinality for ODFCs of full type on $\bbF_q^{2k+1}$ is exactly $q^{k+1}+1$. Consequently, our orbital construction is only two flags away from reaching the mentioned bound.

We devote the rest of the section to complete the construction in Theorem \ref{theo: n=2k+1 ODFC orbita de G} into an ODFC on $\bbF_q^{2k+1}$ with the largest possible size, that is, $q^{k+1}+1$. Notice that  this value coincides with the largest size of a partial spread of $\bbF_q^{2k+1}$ of dimension $k$ or, equivalently, the largest size of a constant dimension code of dimension $k+1$ of $\bbF_q^{2k+1}$ having maximum distance. Therefore, by virtue of Theorem \ref{teo: characterization ODFC a and b}, the problem can be reduced to adding to $\orb_{\bG}(\cF_k)$ and $\orb_{\bG}(\cF_{k+1})$ two appropriate respective subspaces such that the resulting subspace codes of dimensions $k$ and $k+1$ are still of maximum distance. We start tackling this problem for dimension $k$. The following remark will help us in this research. 

\begin{remark}\label{rem: PS Gorla_Ravaganani Orbit + 2 subspaces}
With the notation of (\ref{eq: subspaceU in projected codes n=2k+1}), observe that, in the particular case where $U_1=I_k$ and $U_2=(I_{k+1})_{(k)}$, i.e., the matrix given by the last $k$ rows of $I_{k+1}$, we obtain
\[
 \orb_\bG(\rsp(I_k \ | \ (I_{k+1})_{(k)})) =  \left\lbrace \rsp(I_k \ | \ (M_{k+1}^i)_{(k)}) \ | \ 0\leq i\leq q^{k+1}-2\right\rbrace ,
\]
 where $(M_{k+1}^i)_{(k)}$ is the matrix composed by the last $k$ rows of $M_{k+1}^i$. Therefore, in this case, $\orb_{\bG}(\rsp(I_k \ | \ (I_{k+1})_{(k)}))$ is a subset of the partial spread of $\bbF_q^{2k+1}$ of dimension $k$ given in \cite[Th. 13]{GoRa2014}, which attains the maximum possible size, that is, $q^{k+1}+1$, and can be written as 
\[
\orb_{\bG}(\rsp(I_k \ | \ (I_{k+1})_{(k)})) \cup \{\rsp(I_k \ | \ 0_{k\times (k+1)}),\ \rsp( 0_{k\times k} \ | \ (I_{k+1})_{(k)})\}.
\]
\end{remark}

Inspired by this fact, for any election of full-rank matrices $U_1\in \bbF_q^{k\times k}$ and $U_2\in \bbF_q^{k\times (k+1)}$, we suggest the subspaces
\begin{equation}\label{eq: subspaces U' and U''}
\cU'= \rsp( U_1 \ | \ 0_{k\times (k+1)}) \  \ \text{and} \ \ \cU''= \rsp( 0_{k\times k} \ | \ U_2 ).
\end{equation}
as candidates to make $\orb_\bG(\rsp(U_1 \ | \ U_2))\cup\{\cU', \cU''\}$ be a partial spread of $\bbF_q^{2k+1}$ of dimension $k$. The proof of the following result appears in the final Appendix of the paper.

\begin{proposition}\label{prop: max dist and size dim k n=2k+1}
Let $U_1\in\bbF_q^{k\times k}$ and $U_2\in\bbF_q^{k\times (k+1)}$ matrices such that $\rk(U_1)=\rk(U_2)=k$ and form the $k$-dimensional subspaces $\cU=\rsp(U_1 \ | \ U_2)$ and  $\cU',\cU''$ as in (\ref{eq: subspaces U' and U''}). Then the code $\orb_\bG(\cU)\cup\{\cU', \cU''\}$ is a  partial spread of $\bbF_q^{2k+1}$ of dimension $k$ with cardinality $q^{k+1}+1$.
\end{proposition}

Now we address the same problem for dimension $k+1$. To do so, we consider full-rank matrices $U_1\in\bbF_q^{k\times k}, U_2\in\bbF_q^{k\times (k+1)}$ and vectors $\bv_1\in\bbF_q^k$, $\bv_2\in\bbF_q^{k+1}$ such that the matrix $V=(V_1 \ | \ V_2)$ defined in (\ref{eq: generator matrix V n=2k+1}) has $\rk(V)=\rk(V_2)=k+1$ and put $\cV=\rsp(V)\in \cG_q(k+1, 2k+1)$. In these conditions, by Proposition \ref{prop: n=2k+1 projected code k+1}, the code $\orb_\bG(\cV)$ attains the maximum possible distance and has size $q^{k+1}-1$. Hence, we wonder if it is possible to find two suitable subspaces $\cV'$ and $\cV''$ such that the code $\orb_\bG(\cV)\cup\{\cV', \cV''\}$ still has the maximum distance and achieve the best cardinality. Moreover, in order to use such a code, together with the one given in Proposition \ref{prop: max dist and size dim k n=2k+1}, to construct ODFCs of full type vector on $\bbF_q^{2k+1}$, we also require the condition $\cU'\subset \cV'$ and $\cU''\subset \cV''$. Taking into account the form of the matrix $V$ defined in (\ref{eq: generator matrix V n=2k+1}), it seems quite natural to use subspaces $\cV'$ and $\cV''$, spanned by the rows of matrices
\begin{equation}\label{eq: subspaces V' and V''}
V'=
\begin{pmatrix}
 \begin{array}{c|c}
 U_1 & 0_{k\times (k+1)} \\ \hline
 \bv_1 & \bv_2
 \end{array}
 \end{pmatrix}
\ \text{and} \
V''=
\begin{pmatrix}
 \begin{array}{c|c}
 0_{k\times k}  & U_2  \\ \hline
 \bv_1 & \bv_2
 \end{array}
  \end{pmatrix},
\end{equation}
respectively. Observe that the vector space spanned by the first $k$ rows of $V'$ (resp. $V''$) is precisely $\cU'$ (resp. $\cU''$). The next result states that this pair of subspaces works if, and only if, $\bv_1=\bo_k$. The corresponding proof is also postponed to the final Appendix of the paper.

\begin{proposition}\label{prop: max dist and size dim k+1 n=2k+1}
Let $V$ be  the matrix defined in (\ref{eq: generator matrix V n=2k+1}), taking $\rk(U_1)=\rk(U_2)=k$ and $\bv_2\notin\rsp(U_2)$. Consider $\cV=\rsp(V)$ and subspaces $\cV'$ and $\cV''$ as in (\ref{eq: subspaces V' and V''}). Then the code $\orb_\bG(\cV)\cup\{\cV', \cV''\}$ has the maximum possible distance (i.e., $2k$) if, and only if, $\bv_1=\bo_k$. In such a case, the code attains the largest size, that is, $q^{k+1}+1$.
\end{proposition}

Now, making use of Propositions \ref{prop: max dist and size dim k n=2k+1} and \ref{prop: max dist and size dim k+1 n=2k+1}, we are ready to present the next construction of ODFC of full type on $\bbF_q^{2k+1}$ having the maximum size. 

Take  matrices $U_1\in\bbF_q^{k\times k}$, $U_2\in\bbF_q^{k\times (k+1)}$ such that   $\rk(U_1)=\rk(U_2)=k$ and vectors $\bv_1\in\bbF_q^k$, $\bv_2\in\bbF_q^{k+1}\setminus\rsp(U_2)$ and form subspaces 
\begin{equation}\label{eq: subspaces ODFC max size}
\begin{array}{ll}
 \cU= \rsp(U_1 \ | \ U_2), &  \cV= \rsp\begin{pmatrix}
 \begin{array}{c|c}
 U_1 & U_2 \\ \hline
 \bv_1 & \bv_2
 \end{array}
 \end{pmatrix},\\
 \cU'= \rsp( U_1 \ |\ 0_{k\times(k+1)}), & \cV'= \rsp\begin{pmatrix}
 \begin{array}{c|c}
 U_1 &  0_{k\times(k+1)} \\ \hline
 \bv_1 & \bv_2
 \end{array}
 \end{pmatrix},\\
 \cU''= \rsp( 0_{k\times k} \ | \ U_2), & \cV''= \rsp\begin{pmatrix}
 \begin{array}{c|c}
 0_{k\times k} & U_2 \\ \hline
 \bv_1 & \bv_2
 \end{array}
 \end{pmatrix}.
\end{array}
\end{equation}
With this notation, the next result holds.

\begin{theorem} 
Let $\cF, \cF', \cF''$ be full flags on $\bbF_q^{2k+1}$ such that 
\[
\begin{array}{lll}
\cF_{k}=\cU,    & \cF_{k+1}=\cV,\\
\cF'_{k}=\cU',  & \cF'_{k+1}=\cV',\\
\cF''_{k}=\cU'', &  \cF''_{k+1}=\cV''.
\end{array}
\]
defined as in (\ref{eq: subspaces ODFC max size}). Then the flag code $\cC=\orb_\bG(\cF)\cup\{\cF', \cF''\}$ is an ODFC with the maximum possible cardinality, i.e., $q^{k+1}+1$, if, and only if, $\bv_1=\bo_k$.
\end{theorem}
\begin{proof}
By means of Theorem \ref{teo: characterization ODFC a and b}, we just need to check that the projected codes of dimensions $k$ and $k+1$ attain the maximum distance, which is $2k$ in both cases, and $|\cC|=|\cC_k|=|\cC_{k+1}|$. As proved in Propositions \ref{prop: n=2k+1 projected code k} and \ref{prop: n=2k+1 projected code k+1}, the orbits $\orb_\bG(\cF_k)$ and $\orb_\bG(\cF_{k+1})$ give us the maximum distance. Moreover, by using Proposition \ref{prop: max dist and size dim k n=2k+1}, we obtain that $\cC_k$ attains the maximum distance and size for every choice of $\bv_1$. On the other hand, for $\cC_{k+1}$ this happens if, and only if, $\bv_1=\bo_k$, by Proposition \ref{prop: max dist and size dim k+1 n=2k+1}. In this case, it follows that
\[
|\cC|=|\cC_k|=|\cC_{k+1}|= q^{k+1}+1,
\]
as stated.
\end{proof}


\section{Conclusions and open problems}

In this paper we have addressed the study of flag codes having maximum distance (ODFCs). We have obtained a characterization of such codes in terms of, at most, two of their projected codes.  We have done this first in a general context and then in an orbital scenario (Section \ref{sec:flags}). In particular, these results improve on those obtained in \cite{CasoPlanar} in this respect. 
Next, we have focused on the construction of ODFCs with an orbital structure. To do this, we have used the action of suitable Singer groups. We have provided two different systematic constructions, both of them reaching the maximum possible cardinality. For the first one,  we have exploited the good relationship between Singer groups and Desarguesian spreads to obtain ODFCs having a specific Desarguesian spread as a projected code (Section \ref{subsect:$n=ks$}). For the second construction, we have used the transitive action of Singer groups on hyperplanes and worked with flags of full type vector, thus covering the cases that cannot be considered in the previous construction (Section \ref{subsect:$n=2t+1$}). 

Given that the theoretical results obtained in Section  \ref{sec:flags} do not in any case require working with Desarguesian spreads, a possible research to be done along these lines could include the specific construction of orbital ODFCs having a non-Desarguesian spread among their projected codes. In a wider context, it would be interesting to address the study of flag codes associated with a fixed distance, as well as to provide systematic constructions of them.


\section{Appendix}

\noindent{\em Proof of Proposition \ref{prop: n=2k+1 projected code k}:}
Assume that the code $\orb_\bG(\cU)=\{\cU\cdot  g^i\ |\ 0\leq i\leq q^{k+1}-2\}$ is a partial spread code of dimension $k$ of $\bbF_q^{2k+1}$. In other words, $d_S(\orb_\bG(\cU))=2k$ and $\stab_\bG(\cU)$ is a proper subgroup of $\bG$. In particular, for every $g^j\in\bG\setminus\stab_\bG(\cU)$, it holds that $d_S(\cU, \cU \cdot g^j)=2k$ or, equivalently, $\cU\cap\ \cU\cdot  g^j=\{\bo\}$. Let us see that this necessarily implies that $\rk(U_1)=\rk(U_2)=k$. First, suppose that $\rk(U_2)< k$. Then at least one of the rows of $U_2$ is a linear combination of the other ones. Hence, there must exist a nonzero vector $\ba\in\bbF_q^k$ such that $\ba U_2 = \bo_{k+1}$. Then, for every $0\leq i\leq q^{k+1}-2$, one has that  $\ba U_2 M_{k+1}^i = \bo_{k+1}$ and
\[
\bx= (\ba U_1 \ | \ \bo_{k+1})  = \ba (U_1 \ | \ U_2) = \ba (U_1  \ | \ U_2 M_{k+1}^i)\in \cU\cap\cU\cdot  g^i.
\]
Moreover, since the $k$ rows of $U=(U_1 \ | \ U_2)$ are linearly independent and $\ba\neq\bo_k$, it follows that $ \bx\neq\bo_{2k+1}$. That is a contradiction with $\cU\cap \cU \cdot g^j=\{\bo\}$, for those $g^j\in\bG\setminus\stab_\bG(\cU)$. Hence, it must hold that $\rk(U_2)=k$. Now, let us prove that $U_1$ has also rank $k$. If not, as before, there must exist a nonzero vector $\ba\in\bbF_q^k$ such that $\ba U_1=\bo_k$. Denote $\bx_2=\ba U_2$, which is  a nonzero vector of $\cU_2=\rsp(U_2)$, since the rows of $U_2$ are linearly independent. Besides, notice that $\cU_2$ is a hyperplane of $\bbF_q^{k+1}$. Consider now, another hyperplane (different from $\cU_2$) containing $\bx_2$ too. Since the Singer subgroup $\langle M_{k+1}\rangle$ acts transitively on $\cG_q(k, k+1)$, we can write such a hyperplane as $\cU_2 \cdot M_{k+1}^i,$ for some not scalar matrix $M_{k+1}^i$. In particular, we have $0\neq M_{k+1}^i+I_{k+1} = M_{k+1}^j$, for some $j\in\{ 0, \dots, q^{k+1}-2 \}$.
This, in turn, implies that $M_{k+1}^j$ neither is a scalar matrix. Thus, by Theorem \ref{teo:stabSinger_rectes}, necessarily $\cU_2\neq \cU_2 \cdot M_{k+1}^j$. 
Since $\bo_{k+1}\neq \bx_2\in\cU_2 \cdot M_{k+1}^i = \rsp(U_2M_{k+1}^i)$, there must exist a nonzero vector $\bb\in\bbF_q^k$ such that
\[
\bx_2= \bb U_2M_{k+1}^i =  \bb U_2(M_{k+1}^i +I_{k+1} - I_{k+1}) =  \bb U_2M_{k+1}^j - \bb U_2.
\]
Now, consider the vector
\[
\begin{array}{ccl}
(\ba+\bb)(U_1 \ | \ U_2) & = & (\ba U_1 + \bb U_1 \ | \ \ba U_2 + \bb U_2 )\\
						 & = & ( \bb U_1 \ | \ \bx_2 + \bb U_2M_{k+1}^j - \bx_2) = \bb (U_1 \ | \ U_2M_{k+1}^j).
\end{array}
\]
Finally, since $\bb\neq \bo_k$ and the rows of $(U_1 \ | \ U_2M_{k+1}^j)$ are linearly independent, we have found a nonzero vector lying on $\cU\cap\ \cU\cdot  g^j$. Moreover, since $\cU_2\neq \cU_2 \cdot M_{k+1}^j$, we have that $\cU\neq \cU \cdot g^j$. This represents a contradiction with the fact that $\orb_\bG(\cU)$ is a partial spread. As a result, we conclude that $\rk(U_1)=k$.

Conversely, assume that $\rk(U_1)=\rk(U_2)=k$. Let us see that $\orb_\bG(\cU)$ is a partial spread with $q^{k+1}-1$ elements. To do so, given $g^i\in\bG\setminus \{ I_{2k+1}\}$, consider a vector $\bx$ in $\cU\cap\cU \cdot g^i$. Hence, we can find vectors $\ba, \bb\in\bbF_q^k$ such that 
$$
\bx=\ba(U_1 \ | \ U_2)  = \bb (U_1 \ | \ U_2 M_{k+1}^i).
$$
In particular, it holds $\ba U_1 = \bb U_1$ or, equivalently, $(\ba-\bb)U_1=\bo_k$. Since the matrix $U_1$ is invertible, we conclude that $\ba=\bb$. In this case, we have that $\ba U_2 =\ba U_2 M_{k+1}^i$, i.e., $\ba U_2 (M_{k+1}^i-I_{k+1})= \bo_{k+1}$. Notice that, since $g^i\neq I_{2k+1}$, the matrix $M_{k+1}^i-I_{k+1}$ is invertible and then, it must hold $\ba U_2= \bo_{k+1}$. Moreover, as the $k$ rows of $U_2$ are linearly independent, it follows $\ba=\bo_k$. In other words, the intersection subspace $\cU\cap\cU\cdot  g^i$ is trivial and then $d_S(\cU, \cU\cdot  g^i) = 2k$, for every $g^i\neq I_{2k+1}$. As a consequence $\stab_\bG(\cU)=\{ I_{2k+1} \}$ and $\orb_{\bG}(\cU)$ is a partial spread of cardinality $q^{k+1}-1$.
\cqd

\bigskip

\noindent{\em Proof of Proposition \ref{prop: n=2k+1 projected code k+1}:}
Assume that the code $\orb_\bG(\cV)$ attains the maximum possible distance, i.e., $2k$. Hence,  $\stab_\bG(\cV)$ is a proper subgroup of $\bG$ and, for every $g^i\in\bG\setminus\stab_\bG(\cV)$, it holds $\dim(\cV\cap\cV\cdot  g^i)= 1$. 

Let us start proving that $V_2$ must be an invertible matrix. To do so, arguing by contradiction, we assume that $\rk(V_2)<k+1$. We proceed in two steps:

\begin{enumerate}
\item If $\rk(V_2)\leq k-1 = (k+1)-2$, then we can find, at least, two independent vectors $\ba,\bb\in\bbF_q^{k+1}$ such that $\bb V_2=\bo_{k+1}=\ba V_2$. In this case, both vectors 
\[
\begin{array}{cccccccc}
\bx &=& \ba (V_1 \ | \ V_2 ) &= & \ba (V_1 \ | \ \bo_{k+1} )& = &\ba ( V_1 \ | \ V_2 M_{k+1}^i) & \text{and}\\
\by &=& \bb (V_1 \ | \ V_2 ) &= & \bb (V_1 \ | \ \bo_{k+1} )& = &\bb ( V_1 \ | \ V_2 M_{k+1}^i) &
\end{array}
\]
lie on every subspace of $\orb_\bG(\cV)$. In particular, $\bx, \by\in \cV\cap\cV \cdot g^i$, for every $g^i\in\bG\setminus\stab_\bG(\cV)$. Since in this case we know that $\dim(\cV\cap\cV \cdot g^i)= 1$, there must happen that $\bx=\lambda \by$ for some $\lambda\in\bbF_q$. Then, it holds
\[
\bo_{2k+1} = \bx - \lambda \by = (\ba - \lambda \bb) (V_1 \ | \ V_2).
\]
Moreover, since the rows of $(V_1 \ | \ V_2)$ are linearly independent, we conclude that $\ba=\lambda \bb$, which is a contradiction with the independence of $\ba$ and $\bb$. Hence, it must hold  $k\leq \rk(V_2)\leq k+1$.
\item If $\rk(V_2)=k$, then the subspace $\cV_2=\rsp(V_2)$ is a hyperplane of $\bbF_q^{k+1}$. On the other hand, since $V_1\in\bbF_q^{(k+1)\times k}$, there must exist a nonzero vector $\ba\in\bbF_q^{k+1}$ such that $\ba V_1= \bo_k$. Moreover, notice that the vector $\ba(V_1\ | \ V_2)$ is not zero since $\rk(V_1\ | \ V_2)=k+1$. Hence, $\bx_2=\ba V_2$ is a nonzero vector in the hyperplane $\cV_2$ of $\bbF_q^{k+1}$. Let us consider a different hyperplane of $\bbF_q^{k+1}$ containing $\bx_2$ as well. As the action of $\langle M_{k+1} \rangle$ is transitive on $\cG_q(k, k+1)$, such a hyperplane is of the form $\cV_2 \cdot M_{k+1}^i$, for some not scalar matrix $M_{k+1}^i$. Then we can write $\bx_2 =\bb V_2 M_{k+1}^i$ for some vector $\bb\in\bbF_q^{k+1}$. Observe that  $M_{k+1}^i+I_{k+1}$ is again a power of $M_{k+1}$, say $M_{k+1}^i+I_{k+1}=M_{k+1}^j$, which is not a scalar matrix too. Thus, by Theorem \ref{teo:stabSinger_rectes}, $\cV_2\neq \cV_2 \cdot M_{k+1}^j$  and then $\cV\neq \cV\cdot g^j$. Now, notice that $\bx_2= \bb V_2 M_{k+1}^j - \bb V_2 = \ba V_2$ and 
\[
\bx= (\ba+\bb)(V_1\ | \ V_2) = \bb (V_1\ | \ V_2 M_{k+1}^j) \in \cV\cap\cV \cdot g^j
\]
Moreover, since $\bx_2\neq \bo_{k+1}$, it follows that $\bb V_2\neq \bo_{k+1}$ and then $\bb V_2 M_{k+1}^j\neq \bo_{k+1}$. On the other hand, since $\rk(V_2)=k$, there must exist a nonzero vector $\bc\in\bbF_q^{k+1}$ such that $\bc V_2=\bo_{k+1}$. Observe that, given that $\rk(V_1 \ | \ V_2)=k+1$, then $\bc V_1\neq \bo_{k}$. Hence, the nonzero vector
\[
\by = (\bc V_1\ | \bo_{k+1})= \bc (V_1\ | \ V_2) = \bc (V_1\ | \ V_2 M_{k+1}^j)
\]
lies as well on $\cV\cap\cV \cdot g^j$. We conclude that, $\dim(\cV\cap\cV \cdot g^j)\geq 2$, which contradicts the hypothesis of $\orb_\bG(\cV)$ attaining the maximum possible distance. 
\end{enumerate} 
Hence, assuming that $d_S(\orb_\bG(\cV))=2k$, then it must holds $\rk(V_2)=k+1$. We will see now that $V_1$ needs to have rank equal to $k$. To do so, we assume that $\rk(V_1)\leq k-1$. First of all, notice that if $V_1= 0_{(k+1)\times k}$, since $V_2$ is an invertible matrix, one has that $\orb_\bG(\cV)=\{\cV\}$ and its distance is zero. Hence, we can assume that $1\leq \rk(V_1)\leq k-1$. In this situation, there exist at least two independent vectors $\ba, \bb\in\bbF_q^{k+1}$ such that $\ba V_1 =\bb V_1 =\bo_{k}$. Now, since $\rk(V_2)=k+1,$ the rows of every matrix $V_2 M_{k+1}^i$ span the whole space $\bbF_q^{k+1}$. Take $M_{k+1}^i\neq - I_{k+1}$, then the matrix $M_{k+1}^i+I_{k+1}$ is again a power of $M_{k+1}$, say $M_{k+1}^i+I_{k+1}=M_{k+1}^j,$ for some $0\leq j \leq q^{k+1}-2$. Now we express both (nonzero) vectors $\ba V_2$ and $\bb V_2$ as linear combinations of the rows of $V_2M_{k+1}^i$ and obtain the following equalities:
\begin{equation}\label{eq: vectors aV2 and bV2}
\begin{array}{ccccccc}
\ba V_2 &=& \bc V_2 M_{k+1}^i &=&  \bc V_2 (M_{k+1}^i+I_{k+1} - I_{k+1}) &=&  \bc V_2 M_{k+1}^j - \bc V_2,\\
\bb V_2 &=& \bd V_2 M_{k+1}^i &=&  \bd V_2 (M_{k+1}^i+I_{k+1} - I_{k+1}) &=&  \bd V_2 M_{k+1}^j - \bd V_2,
\end{array}
\end{equation}
for some nonzero vectors $\bc, \bd \in\bbF_q^{k+1}$. Hence, both vectors
\[
\begin{array}{ccccccc}
\bx &=& (\ba +\bc)(V_1\ | \ V_2) &=& (\bc V_1 \ | \ \ba V_2 +\bc V_2) &=& \bc (V_1 \ | \  V_2 M_{k+1}^j),\\
\by &=& (\bb +\bd)(V_1\ | \ V_2) &=& (\bd V_1 \ | \ \bb V_2 +\bd V_2) &=& \bd (V_1 \ | \  V_2 M_{k+1}^j),
\end{array}
\]
lie on $\cV\cap\cV\cdot  g^j$. Let us see that they are independent. Otherwise,  $\bx =\lambda \by$, for some $\lambda\in \bbF_q$ and it must hold
\[
\bo_{k+1}= (\bc-\lambda\bd) V_2 M_{k+1}^j.
\]
Since $V_2 M_{k+1}^j$ is invertible, we conclude that $\bc=\lambda\bd$. As a result, by (\ref{eq: vectors aV2 and bV2}), we obtain
\[
\ba V_2 = \bc V_2 M_{k+1}^i = \lambda \bd V_2 M_{k+1}^i = \lambda \bb V_2 
\]
and, given that $\rk(V_2)=k+1$, it follows $\ba=\lambda \bb$, which is not possible. Hence,  $\dim(\cV\cap\cV \cdot g^j)\geq 2$. Moreover, let us see that $\cV\neq \cV \cdot g^j$. Note that $\dim(\cV + \cV \cdot g^j)$ is exactly the rank
\[
\rk
\begin{pmatrix}
\begin{array}{c|c}
V_1 & V_2\\ \hline
V_1 & V_2 M_{k+1}^j 
\end{array}
\end{pmatrix}
=
\rk
\begin{pmatrix}
\begin{array}{c|c}
V_1 & V_2\\ \hline
0_{(k+1)\times k} & V_2 M_{k+1}^i
\end{array}
\end{pmatrix}
\geq
\rk(V_1)+\rk( V_2 M_{k+1}^i) 
\geq k+2,
\]
which is greater than $\dim(\cV)=\dim(\cV\cdot  g^j)=k+1$. Hence, $\cV\neq \cV \cdot g^j$ and we obtain a contradiction with $d_S(\orb_\bG(\cV))=2k$. Then, it must happen $\rk(V_1)=k$ and $\rk(V_2)=k+1.$ 

Now, let us prove that the converse is also true. Assume that $\rk(V_1)=k$ and $\rk(V_2)=k+1$ and take $g^i\in\bG\setminus\{I_{2k+1}\}$. We show that $\dim(\cV\cap \cV\cdot g^i)=1$. To do so, consider two arbitrary vectors $\bx,\by\in\cV\cap\cV \cdot g^i$ and write them as 
$$
\begin{array}{ccccl}
\bx &=& \ba(V_1 \ | \ V_2) &=& \bc(V_1 \ | \ V_2M_{k+1}^i) \ \text{and}\\
\by &=& \bb(V_1 \ | \ V_2) &=& \bd(V_1 \ | \ V_2M_{k+1}^i),  
\end{array}
$$
for some vectors $\ba, \bb, \bc, \bd\in\bbF_q^{k+1}$. Observe that $(\ba-\bc)V_1= (\bb-\bd) V_1=\bo_k.$ Since $V_1\in \bbF_q^{(k+1)\times k}$ with $\rk(V_1)=k$, it follows that $\ba-\bc$ and $\bb-\bd$ must be proportional vectors. Let us write $(\ba-\bc)=\lambda(\bb-\bd),$ for some $\lambda \in\bbF_q$. Equivalently, $\ba-\lambda \bb = \bc-\lambda\bd$.  Moreover, since $
\ba V_2= \bc V_2 M_{k+1}^i$ and $\bb V_2= \bd V_2 M_{k+1}^i$, we obtain
\[
(\ba-\lambda\bb) V_2 = (\bc -\lambda\bd) V_2 M_{k+1}^i =  (\ba-\lambda \bb)V_2 M_{k+1}^i,
\]
or, equivalently, $(\ba-\lambda\bb)V_2(M_{k+1}^i-I_{k+1})=\bo_{k+1}$. Last, since $g^i\neq I_{2k+1}$, then $V_2(M_{k+1}^i-I_{k+1})$ is an invertible matrix and it follows $\ba=\lambda\bb$ and $\bx=\lambda\by$. As a result, for every $g^i\in\bG\setminus\{I_{2k+1}\}$, it occurs
$\dim(\cV\cap\cV \cdot g^i)=1$ or, equivalently, $d_S(\cV, \cV\cdot  g^i) = 2k$. Consequently, $\stab_\bG(\cV)=\{I_{2k+1}\}$ and $\orb_\bG(\cV)$ has maximum distance and cardinality $q^{k+1}-1$.
\cqd

\bigskip

\noindent{\em Proof of Proposition \ref{prop: max dist and size dim k n=2k+1}:}
Since $\rk(U_1)=\rk(U_2)= k$, Proposition \ref{prop: n=2k+1 projected code k} proves that the orbit $\orb_\bG(\cU)$ is a $k$-partial spread of $\bbF_q^{2k+1}$. Hence, it suffices to see that adding the two subspaces $\cU'$ and $\cU''$ defined in (\ref{eq: subspaces U' and U''}) does not decrease the distance. Observe that two $k$-dimensional subspaces in $\bbF_q^{2k+1}$ attain the maximum possible distance if, and only if, they intersect trivially or, equivalently, if their sum subspace has dimension $2k$. It is clear that
$$
\dim(\cU \cdot g^i + \cU')
=
\rk
\begin{pmatrix}
\begin{array}{c|c}
U_1 & U_2 M_{k+1}^i\\ \hline
U_1 & 0_{k\times (k+1)}
\end{array} 
\end{pmatrix}
=
\rk
\begin{pmatrix}
\begin{array}{c|c}
0_{k\times k} & U_2 M_{k+1}^i\\ \hline
U_1 & 0_{k\times (k+1)}
\end{array} 
\end{pmatrix}
= 2k,
$$
$$
2k \geq \dim(\cU \cdot g^i + \cU'')
=
\rk
\begin{pmatrix}
\begin{array}{c|c}
U_1 & U_2 M_{k+1}^i\\ \hline
0_k & U_2
\end{array} 
\end{pmatrix}
\geq
\rk(U_1)+\rk(U_2)= 2k
$$
and
$$
\dim(\cU' + \cU'')
=
\rk
\begin{pmatrix}
\begin{array}{c|c}
U_1 & 0_{k\times (k+1)}\\\hline
0_k & U_2
\end{array} 
\end{pmatrix}
=
\rk(U_1)+\rk(U_2) = 2k.
$$
As a result, in these three situations we obtain distance equal to $2k$ and we conclude that the code $\orb_\bG(\cU)\cup\{\cU', \cU''\}$ is a partial spread with $q^{k+1}+1$ elements.
\cqd

\bigskip

\noindent{\em Proof of Proposition \ref{prop: max dist and size dim k+1 n=2k+1}:}
Notice that, by virtue of Proposition \ref{prop: n=2k+1 projected code k+1}, we have that $d_S(\orb_\bG(\cV))=2k$. Hence, in order to give the minimum distance of the code $\orb_\bG(\cV)\cup\{\cV', \cV''\},$ we just need to compute distances between pairs of different subspaces not both in $\orb_\bG(\cV)$. Moreover, recall that two $(k+1)$-dimensional subspaces of $\bbF_q^{2k+1}$ have the maximum possible distance if, and only if, their sum is the whole vector space $\bbF_q^{2k+1}$.

Let us start proving that taking $\bv_1=\bo_k$ leads to a construction of maximum distance. To do so, just notice that
\[
\begin{array}{ccl}
\dim(\cV\cdot  g^i + \cV')  =  \rk
\begin{pmatrix}
\begin{array}{c|c}
U_1 & U_2M_{k+1}^i \\ \hline
\bo_k & \bv_2 M_{k+1}^i \\ \hline
U_1 & 0_{k\times (k+1)} \\ \hline
\bo_k & \bv_2
\end{array}
\end{pmatrix}
& =&
\rk
\begin{pmatrix}
\begin{array}{c|c}
0_{k\times k} & U_2M_{k+1}^i \\ \hline
\bo_k & \bv_2 M_{k+1}^i \\ \hline
U_1 & 0_{k\times (k+1)}
\end{array}
\end{pmatrix} \\
& = & \rk(U_1)+ \rk(V_2M_{k+1}^i)= 2k+1,
\end{array}
\]

\[
\dim(\cV \cdot g^i + \cV'')=
\rk
\begin{pmatrix}
\begin{array}{c|c}
U_1 & U_2M_{k+1}^i \\ \hline
\bo_k & \bv_2 M_{k+1}^i \\ \hline
0_{(k+1)\times k} & U_2 \\ \hline
\bo_k & \bv_2 
\end{array}
\end{pmatrix}
\geq 
\rk(U_1)+\rk(V_2)=2k+1
\]
and
\[
\begin{array}{ccl}
\dim(\cV'+ \cV'')=
\rk
\begin{pmatrix}
\begin{array}{c|c}
U_1 & 0_{k\times (k+1)} \\ \hline
\bo_k & \bv_2\\ \hline
0_{(k+1)\times k} & U_2 \\ \hline
\bo_k & \bv_2
\end{array}
\end{pmatrix}
&=& \rk(U_1)+\rk(V_2)= 2k+1.
\end{array}
\]
Hence, we obtain that $\cV \cdot g^i + \cV' = \cV \cdot g^i+ \cV''= \cV'+\cV''=\bbF_q^{2k+1},$  for every $0\leq i \leq q^{k+1}-2$.  Consequently, the distance of the code is the maximum possible one. In particular, $\cV'$ and $\cV''$ are different and they do not lie in the orbit of $\cV$. Thus, it follows $|\orb_\bG(\cV)\cup\{\cV', \cV''\}|=q^{k+1}+1$, i.e., the largest size for its distance. 

Conversely, we show that taking $\bv_1=\bo_k$ is the only possibility for $\orb_\bG(\cV)\cup \{ \cV', \cV'' \}$ to attain the maximum distance. To do so, assume that $\bv_1$ is a nonzero vector of $\bbF_q^k$. Hence, since $\rk(U_1)=k$, there exists a nonzero vector $\ba\in\bbF_q^k$ such that $\bv_1=\ba U_1$. We will exhibit a explicit subspace $\cV \cdot g^j\in \orb_{\bG}(\cV)$ such that $d_S(\cV\cdot  g^j,\cV'')<2k$.

Since $U_2\in\bbF_q^{k\times (k+1)}$ has $\rk(U_2)=k$, the subspace $\cU_2=\rsp(U_2)$ is a hyperplane of $\bbF_q^{k+1}$ and, as $\rk(V_2)=k+1$, it follows that $\bv_2\notin \cU_2$. In particular, $\bv_2$ and $\ba U_2$ are linearly independent vectors and so they are $\bv_2$ and $\bv_2-\ba U_2$. Consider a hyperplane of $\bbF_q^{k+1}$ containing the vector $\bv_2-\ba U_2$ but not $\bv_2$. Recall that, since the Singer subgroup generated by $M_{k+1}$ acts transitively on $\cG_q(k, k+1)$ (see Theorem \ref{teo:stabSinger_rectes}), such a hyperplane is of the form $\cU_2 \cdot M_{k+1}^i$, for some $i\in \{0,\dots, q^{k+1}-2\}$. Observe that, for this choice of $i$, we have $\bv_2-\ba U_2\in\cU_2 \cdot M_{k+1}^i$ and $\bv_2\notin\cU_2\cdot  M_{k+1}^i$. Thus, at this point, we have that 
\[
(\bv_2-\ba U_2)M_{k+1}^{-i}\in\cU_2 \quad \mbox{ and } \quad \bv_2,\ \bv_2 M_{k+1}^{-i}\notin\cU_2.
\]

On the other hand, as $V_2\in\GL(k+1,q)$, every vector in $\bbF_q^{k+1}$ can be written as a linear combination of its rows, which are the ones of $U_2$ together with $\bv_2$. In particular, there must exist $\lambda\in\bbF_q$ and $\bb\in\bbF_q^k$ such that  $\bv_2M_{k+1}^{-i} = \lambda \bv_2 + \bb U_2.$ Moreover, since $\bv_2M_{k+1}^{-i}\notin\cU_2$, it follows $\lambda \neq 0$. Then the matrix $\lambda^{-1}M_{k+1}^{-i}\in\bbF_q[M_{k+1}]$ is a power of $M_{k+1}$ and we can write $\lambda^{-1}M_{k+1}^{-i}= M_{k+1}^j$, for certain exponent $j\in  \{0,\dots, q^{k+1}-2\}$.  Observe that, for this matrix $M_{k+1}^j$, we have that $\bv_2M_{k+1}^j- \bv_2\in \cU_2$  and also $(\bv_2-\ba U_2)M_{k+1}^j\in\cU_2$. As a result, their difference, i.e., the vector $\bv_2-\ba U_2 M_{k+1}^j\in\cU_2$ as well.
Now, consider the subspace $\cV \cdot g^j$. Let us see that $d_S(\cV \cdot g^j, \cV'')$ is not the maximum one or, equivalently, that $\cV \cdot g^j+\cV''\neq \bbF_q^{2k+1}$. Observe that 
\[
\begin{array}{ccl}
\dim(\cV \cdot g^j+\cV'') &=& 
\rk
\begin{pmatrix}
\begin{array}{c|c}
U_1 & U_2M_{k+1}^j \\ \hline
\bv_1 & \bv_2 M_{k+1}^j \\ \hline
0_{(k+1)\times k} & U_2 \\ \hline
\bv_1 & \bv_2 
\end{array}
\end{pmatrix} \\
& = & \rk
\begin{pmatrix}
\begin{array}{c|c}
U_1 & U_2M_{k+1}^j \\ \hline
\bo_k & \bv_2 M_{k+1}^j - \ba U_2M_{k+1}^j \\ \hline
0_{(k+1)\times k} & U_2 \\ \hline
\bo_k & \bv_2 - \ba U_2M_{k+1}^j
\end{array}
\end{pmatrix}\\
& = & \rk
\begin{pmatrix}
\begin{array}{c|c}
U_1 & U_2M_{k+1}^j \\ \hline
0_{(k+1)\times k} & U_2 
\end{array}
\end{pmatrix}\\
& = & \dim(\cU\cdot g^j + \cU'')=2k.
\end{array}
\]
Since $k>1$, it follows that $2k>k+1$. Consequently, for this precise value of $j$, it holds that $\cV \cdot g^j\neq \cV''$ and $d_S(\cV \cdot g^j, \cV'')< 2k$. We conclude that $\orb_\bG(\cV)\cup\{\cV', \cV''\}$ does not attain the maximum distance.
\cqd

\end{document}